\documentclass[12pt]{article}
\usepackage{cite}
\usepackage{mathrsfs}
\usepackage{graphicx}
\usepackage{float}
\usepackage{curves}
\usepackage{amsfonts}
\usepackage{amsmath}
\usepackage{amsfonts,amssymb,color}
\usepackage{dsfont}
\usepackage{graphicx}
\usepackage{curves}
\usepackage{mathrsfs}
\usepackage{pifont}
\usepackage{amssymb}

\numberwithin{equation}{section}
\newtheorem{thm}{Theorem}[section]
\newtheorem{prop}{Proposition}[section]
\newtheorem{defn}{Definition}[section]
\newtheorem{lem}{Lemma}[section]
\newtheorem{cor}{Corollary}[section]
\newtheorem{rem}{Remark}[section]

\def\emptyset{\mbox{{\rm \O}}}

\newenvironment {pf} {\noindent {\bf Proof.}}{\hspace*{\fill}$\Box$\par\vspace{3mm}}
\begin{document} 

\baselineskip=30pt

\title{\bf \Large The general spectral radius and majorization theorem of $t$-cone graphs with given degree sequences\footnotemark[1]\\[5mm]}

\author{Yufei Huang$^{1}$, Muhuo Liu$^{2}$~
\\[2mm] \small $^1$ Department of Mathematics Teaching, Guangzhou Civil Aviation College, Guangzhou, 510403
\\\small $^2$ Department of Mathematics, South China Agricultural University, Guangzhou, 510642}

\renewcommand{\thefootnote}{\fnsymbol{footnote}}

\footnotetext[1]{This work is supported by NNSF of China (Nos. 11571123, 11501139),  and  Guangdong Province Ordinary University Characteristic Innovation Project
(No.2017KTSCX020). E-mail address: liumuhuo@163.com (M. Liu, Corresponding author), fayger@qq.com (Y. Huang).}

\date{}
\maketitle
\date{}


\begin{abstract}
The general spectral radius of a graph $G$, denoted by $\Theta(G,\alpha)$,
is the maximal eigenvalue of $M_{\alpha}(G)=A(G)+\alpha D(G)$ $(\alpha\geq 0)$,
where $A(G)$ and $D(G)$ are the adjacency matrix and the diagonal matrix of vertex degrees of $G$, respectively.
A graph $G$ is called $\Theta_\alpha$-maximal
in a class of connected simple graphs $\mathcal {G}$
if $\Theta(G,\alpha)$ is maximal among all graphs of $\mathcal {G}$.
A $t$-cone $c$-cyclic graph is the join of a complete graph $K_t$ and a $c$-cyclic connected simple graph.
Let $\pi=\big(d_1,d_2,\ldots,d_n\big)$ and $\pi'=\big(d'_1,d'_2,\ldots,d'_n\big)$ be two non-increasing degree sequences of $t$-cone $c$-cyclic graphs with $n$ vertices.
We say $\pi$ is strictly  majorized by $\pi'$, denoted by $\pi \lhd \pi'$,
if $\pi\neq \pi'$, $\sum_{i=1}^n d_i=\sum_{i=1}^n
d_i'$, and $\sum_{i=1}^k d_i\leq \sum_{i=1}^k
d_i'$ for $k=1,2,\ldots,n-1$.
Denote by $\Gamma(\pi,t;c)$ the class of $t$-cone $c$-cyclic graphs with $\pi$ as its degree sequence.
In this paper,
we determine some properties of $\Theta_\alpha$-maximal graphs of $\Gamma(\pi,t;c)$ and characterize the unique $\Theta_\alpha$-maximal graph of $\Gamma(\pi,t;0)$ \big(resp. $\Gamma(\pi,t;1)$ and $\Gamma(\pi,t;2)$\big).
Moreover,
we prove that if $\pi \lhd \pi'$, $G$ and $G'$ are the $\Theta_\alpha$-maximal graphs of
$\Gamma(\pi,t;c)$ and $\Gamma(\pi',t;c)$  respectively,
then $\Theta(G,\alpha)<\Theta(G',\alpha)$ for $c\in \big\{0,1\big\}$, and we also consider the  similar result for $c=2$. \\



\end{abstract}

\baselineskip=0.3in



\section{\large{Introduction}}

Let
$G=\big(V(G),E(G)\big)$
be a connected simple graph with $n$ vertices and $m$ edges.
If $m=n+c-1$, then
$G$ is called a $c$-{\em cyclic graph}, where $c$ is an
integer with $0\leq c\leq {n \choose 2}-n+1$.
Particularly, if $c=0,1$ and $2$, then $G$ is
called a {\em tree}, {\em unicyclic graph} and {\em bicyclic graph}, respectively. Let  $N_G(u)$ and $d_G(u)$, short for $N(u)$ and $d(u)$ if there is no risk of confusion,  denote the neighbor set and degree  of vertex $u$ of $G$, namely, $d(u)=|N(u)|$.
Denote by $\delta(G)$ the minimum degree of $G$.
Hereafter, we suppose that $V(G)=\big\{v_1,v_2,\ldots,v_n\big\}$.
If $d_i=d(v_i)$ for $1\leq i\leq n$,
then we call the sequence $\pi=\left(d_1,d_2,\ldots,d_n\right)$ the {\em degree sequence} of $G$.
Throughout this paper,
we enumerate the degrees in non-increasing order,
i.e. $d_1\geq d_2\geq\cdots \geq d_n$.
Consequently,
$d(v_1)\geq d(v_2)\geq\cdots \geq d(v_n)$.
Hereafter,
let $\Gamma(\pi)$ be the class of connected simple graphs
with $\pi$ as its degree sequence.
Let $C_n$, $K_n$ and $P_n$ denote the cycle,
complete graph and path with $n$ vertices, respectively.
For all undefined  notations and terminologies of this paper,
the readers may refer to
\cite{Bondy,Rowlinson}.

As usual, we use  $A(G)$ and $D(G)$ to denote the adjacency matrix and
the diagonal matrix of vertex degrees of $G$, respectively. Then,
the matrix $Q(G)=A(G)+D(G)$ is called the {\em signless Laplacian
matrix} of $G$. As $A(G)$ and $Q(G)$ are both real symmetric, we can use
$\rho(G)$ and $\mu(G)$, respectively, to denote  the {\em spectral radius} and {\em signless Laplacian
spectral radius} of $G$, that is, the maximum eigenvalues of $A(G)$ and $Q(G)$,
respectively.  As a natural extension to $A(G)$ and $Q(G)$, Liu et al. \cite{LiuQ2014} constructed  the {\em general matrix} $M_\alpha(G)=A(G)+\alpha D(G)$   for a graph $G$.
In view of the real symmetry of  $M_{\alpha}(G)$,
we use $\Theta(G,\alpha)$
to denote the maximal  eigenvalue  of $M_{\alpha}(G)$,
and call  $\Theta(G,\alpha)$ the {\em general spectral radius} of $G$.
For convenience,
if $\Theta(G,\alpha)$ is maximal among all graphs in a class of connected simple graphs $\mathcal {G}$,
then we call $G$ as a  $\Theta_\alpha$-{\em maximal graph} of $\mathcal {G}$.
From these involving definitions,  one can easily see that $M_{0}(G)=A(G)$ and $M_{1}(G)=Q(G)$,
and so
$\Theta(G,0)=\rho(G)$ and  $\Theta(G,1)=\mu(G)$. With the similar reason,
$\Theta(G,0)$-maximal graph and $\Theta(G,1)$-maximal graph are always simplify  called $\rho$-{\em maximal graph} and $\mu$-{\em maximal graph},
respectively.

Let
{\bf y}$=\big(y_1,y_2,\ldots,y_n\big)$ and {\bf z}$=\big(z_1,z_2,\ldots,z_n\big)$ be two non-increasing sequences of real numbers. If
$$
\sum_{i=1}^k y_i\leq \sum_{i=1}^k
z_i~~\text{for}~~k=1,2,\ldots,n-1,
~~\text{and}~~\sum_{i=1}^n y_i=\sum_{i=1}^n z_i,
$$
then {\bf y} is said to be {\em majorized} by {\bf z} \cite{Marshall},
and denoted by {\bf y}$\unlhd $ {\bf z}.
If {\bf y}$\unlhd $ {\bf z} and {\bf y}$\neq$ {\bf z},
then {\bf y}  is called  {\em strictly majorized} by {\bf z},
which will be denoted by {\bf y}$\lhd $ {\bf z} hereafter.

Extremal results on spectral radius and/or signless Laplacian spectral radius in some fixed graph categories have received much attention in recent years.
In this line,
B{\i}y{\i}ko\u{g}lu and Leydold \cite{Biyikoglu} firstly considered the extremal problem of spectral radius in a class of connected simple graphs with given degree sequence, and   they also proved  the majorization theorem between two $\rho$-maximal trees with different tree degree sequences, that is,

\begin{thm}{\em\cite{Biyikoglu}}
Let $\pi$ and $\pi$ be two different non-increasing degree sequences of trees
with $\pi\lhd \pi'$. If $T$ and $T'$ are, respectively,  the $\rho$-maximal trees of  $\Gamma(\pi)$ and
$\Gamma(\pi')$, then  $\rho(T)<\rho(T')$.
\end{thm}

Simultaneously,
Zhang \cite{Zhang0} proved the similar  majorization theorem for the (signless) Laplacian spectral radius of trees,
and characterized the corresponding unique $\mu$-maximal tree of $\Gamma(\pi)$ for any prescribed   tree degree sequence $\pi$.
Since then,
similar problems have been studied extensively.
The majorization theorems for
(signless Laplacian) spectral radius
of unicyclic graphs and bicyclic graphs were also discovered,
and the unique $\rho$-maximal ($\mu$-maximal) graph of $\Gamma(\pi)$ were characterized for
unicyclic graph and bicyclic graph degree sequences $\pi$,
respectively \cite{Liu2,Zhang1,Jiang,Huang,LiuY2010,FB2010}.
Here we have to point out that the majorization theorem of (signless Laplacian) spectral radius can not hold for all $c$-cyclic graphs, as  counterexamples show that the majorization theorem for the (resp., signless Laplacian) spectral radius of 3-cyclic  (resp., 4-cyclic) graphs does not hold  \cite{Jiang,Liu7}.
Besides,
Liu et al. \cite{Liu2,Liu4,Liu5} proved
the majorization theorems for (signless Laplacian) spectral radius of $c$-cyclic graphs with additional restrictions
and they also proved the majorization theorem for (signless Laplacian) spectral radius of pseudographs \cite{Liu6}.
For more results on this field, one may
refer to \cite{Liu1,Liu7,Zhangx}.

Let $G_1$ and $G_2$ be two vertex-disjointed graphs.
Let $G_1\vee G_2$ denote the {\em join graph} of $G_1$ and $G_2$,
which is obtained by joining each vertex of $G_1$ and each vertex of $G_2$.
Let $n$ and $t$ be two integers with $0\leq t\leq n-2$, and let $c$   be an integer with $0\leq c\leq {n-t \choose 2}-n+t+1$. If $G$ is a connected simple ($c$-cyclic)  graph with $n-t$ vertices,
then $K_t\vee G$ is called a $t$-{\em cone} ($c$-{\em cyclic$)$ graph} with $n$ vertices.
Here,
$K_0$ denotes the null graph,
which is defined as a graph with empty vertex set.
According to the above definitions,
$K_0\vee G=G$ which is the usual ($c$-cyclic) graph,
and
$K_1\vee G$ is also  called a {\em single-cone} ($c$-{\em cyclic$)$ graph} in \cite{Luo}.
For simplification,
$t$-cone $0$-cyclic graph, $t$-cone $1$-cyclic graph and
$t$-cone $2$-cyclic graph,
are also called the $t$-{\em cone tree}, $t$-{\em cone unicyclic graph} and
$t$-{\em cone bicyclic graph}, respectively.

Let $\pi$ be a degree sequence of $t$-cone $c$-cyclic graph.
Different from the definition of $\Gamma(\pi)$,
denote by $\Gamma(\pi,t;c)$ the class of $t$-cone $c$-cyclic graphs with $\pi$ as its degree sequence.
Here we would like to point out  that $K_t\vee H\in \Gamma(\pi)$
can not guarantee the connectivity of $H$ and hence $\Gamma(\pi,t;c)\subseteq \Gamma(\pi)$.
For simplification, we
denote by $\Gamma(\pi,t;0)=\mathscr{T}(\pi,t)$,
$\Gamma(\pi,t;1)=\mathscr{U}(\pi,t)$
and $\Gamma(\pi,t;2)=\mathscr{B}(\pi,t)$.

Recently,
Luo et al. \cite{Luo} proved the majorization theorems for the (signless Laplacian) spectral radius of single-cone trees and single-cone unicyclic graphs respectively,
but they characterize neither the $\rho$-maximal graphs nor the $\mu$-maximal graphs of $\mathscr{T}(\pi,1)$ and $\mathscr{U}(\pi,1)$.
Motivated by their research but further than that,
in this paper,
we first characterize the $\Theta_\alpha$-maximal graphs of $\Gamma(\pi,t;c)$,
and then consider the majorization theorems for the general spectral radius of $t$-cone trees, $t$-cone unicyclic graphs and $t$-cone bicyclic graphs,
respectively.

The rest of this paper is organized as follows.
In Section 2,
we recall some basic notions and lemmas which are useful in the proof of our main results.
In Section 3,
some properties of $\Theta_\alpha$-maximal graphs of $\Gamma(\pi,t;c)$ are determined,
and they will play an important role in the characterization of
the unique  $\Theta_\alpha$-maximal graph  of $\mathscr{T}(\pi,t)$, $\mathscr{U}(\pi,t)$ and/or $\mathscr{B}(\pi,t)$, respectively,
in Section 4.
In Section 5,
the majorization theorems for the general spectral radius of $t$-cone trees,
$t$-cone unicyclic graphs and $t$-cone bicyclic graphs are considered,
respectively.
Finally, we present many    related published  results in Section 6,
which can be deduced from our new results straightly.

In what follows \big(including all results\big),
{\bf for the sake of simplification, unless special indicated},
$\alpha$ always defines a non-negative real number,
$t$ and $n$ are two integers  such that $0\leq t\leq n-2$;
let $\pi=\big(d_1,d_2,\ldots,d_n\big)$ and $\pi'=\big(d'_1,d'_2,\ldots,d'_n\big)$
be two non-increasing degree sequences of $t$-cone $c$-cyclic graphs such that $\pi\lhd \pi'$,
and denote by $G_{\pi}$ and $G_{\pi'}$ the $\Theta_\alpha$-maximal graphs of $\Gamma(\pi,t;c)$ and $\Gamma(\pi',t;c)$, respectively;
denote by $\pi^*=\left(d^*_{t+1},d^*_{t+2},\ldots,d^*_{n}\right)$, where  $d^*_j=d_j-t$ for $j\in \big\{t+1,t+2,\ldots,n\big\}$,
then $\pi^*$ is the degree sequence of a $c$-cyclic graph with $n-t$ vertices;
besides, by $G=K_t\vee H\in \Gamma(\pi,t;c)$ or $G=K_t\vee H$ being a $t$-cone graph with $n$ vertices,
we means that $H$ is a $c$-cyclic graph with $n-t$ vertices,
and denoted by $V(H)=\big\{v_{t+1},v_{t+2},\ldots,v_n\big\}$;
we use $P_G(u,v)$ to denote a shortest path connecting the vertices $u$ and $v$ in $G$,
and  $dist_G(u,~v)$ to denote the distance between $u$ and $v$ in $G$,
namely, $dist_G(u,v)$ is equal to the number of edges of $P_G(u,v)$.


\section{\large{Preliminaries}}

Let $G$ be a connected simple graph with $n$ vertices.
Let ${\bf \varphi}=\big(\varphi(v_1),\varphi(v_2),\ldots,\varphi(v_n)\big)^T\in R^n$ be a unit column vector defined on $V(G)=\big\{v_1,v_2,\ldots,v_n\big\}$.
Then, the {\em Rayleigh quotient} \cite{Horn} of the general matrix $M_{\alpha}(G)=A(G)+\alpha D(G)$ is defined as:
$$
\mathcal{R}_{M_{\alpha} (G)}({\bf \varphi})={\bf \varphi}^T M_{\alpha}(G) {\bf \varphi}=2\sum_{uv\in E(G)} \varphi(u)\varphi(v)+\alpha \sum_{v\in V(G)} d(v)\,\varphi^2(v).
$$
When $G$ is connected, since  $\alpha  \geq 0$ and $M_{\alpha}(G)$ is a  nonnegative irreducible matrix,
by the famous  Perron-Frobenius theorem for nonnegative irreducible matrix \cite{Horn},
there exists a unique unit positive eigenvector ${\bf f}=\big(f(v_1),f(v_2),\ldots,f(v_n)\big)^T$
such that $M_{\alpha}(G) {\bf f}=\Theta(G,\alpha) {\bf f}$,
and
\begin{equation}\label{21e}
\Theta(H,\alpha)<\Theta(G,\alpha)~~\text{holds~for any proper subgraph}~~H\subset G.
\end{equation}
Moreover, this eigenvector ${\bf f}$ is called the {\em Perron vector} of $M_{\alpha}(G)$.
If not specified in the following,
we always use the notation ${\bf f_G}=\big(f_G(v_1),f_G(v_2),\ldots,f_G(v_n)\big)^T$ to denote the Perron vector of $M_{\alpha}(G)$,\underline{}
and $f_G(v)$ is called the {\em $\Theta_{\alpha}$-weight} of vertex $v$.
By the Rayleigh-Ritz theorem \cite{Horn},
for any unit column vector ${\bf \varphi}$ defined on $V(G)$, we have
$$\Theta(G,\alpha)\geq \mathcal{R}_{M_{\alpha}(G)}({\bf \varphi}),$$
where the equality holds if and only if ${\bf \varphi}={\bf  f_G}$, namely, ${\bf \varphi}$ is the Perron vector of $G$.
Moreover,
for every vertex $v\in V(G)$, we have
\begin{equation}\label{22e}
\Theta(G,\alpha) f_G(v)=\big(M_{\alpha}(G) {\bf f_G}\big)(v)=\sum_{u\in N(v)} f_G(u)
+\alpha  d_G(v)f_G(v).
\end{equation}

The following two theorems on graph operations, usually called shifting and switching respectively, are very useful in the research of extremal theory of graph spectrum.

\begin{thm}\label{21t}{\em\cite{Liu1}}
Let $u,v$ be two vertices of a connected simple graph $G$,
and $w_1,w_2,\ldots,w_k$ $\big(1\leq k\leq d(v)\big)$
be some vertices of $N(v)\setminus \big(N(u)\cup \{u\}\big)$.
Let $G'=G+w_1u+w_2u+\cdots+w_ku-w_1v-w_2v-\cdots-w_kv$.
If $f_G(u)\geq f_G(v)$,
then $\Theta(G',\alpha)>\Theta(G,\alpha)$.
\end{thm}

\begin{thm}\label{22t}{\em\cite{Liu1}}
Let $G$ be a connected simple graph such that
$uv\in E(G)$,
$xy\in E(G)$,
$uy\notin E(G)$ and $xv\notin E(G)$.
Let $G'=G+uy+xv-uv-xy$.
If $G'$ is connected,
$f_G(u)\geq f_G(x)$ and
$f_G(y)\geq f_G(v)$,
then $\Theta(G',\alpha)\geq\Theta(G,\alpha)$,
with equality if and only if $f_G(u)=f_G(x)$ and $f_G(y)=f_G(v)$.
\end{thm}

\begin{cor}\label{23c}
Let $G=K_t\vee H$ be a $t$-cone graph.
Let $u,v\in V(H)$ such that  $w_1,w_2,\ldots,w_k$ $\big(1\leq k\leq d_H(v)\big)$
are  $k$ vertices of $N_H(v)\setminus \big(N_H(u)\cup \{u\}\big)$.
Let $H'=H+w_1u+w_2u+\cdots+w_ku-w_1v-w_2v-\cdots-w_kv$ and $G'=K_t\vee H'$.
If $f_G(u)\geq f_G(v)$,
then $\Theta(G',\alpha)>\Theta(G,\alpha)$.
\end{cor}

\begin{pf}
Obviously,
$G'=G+w_1u+w_2u+\cdots+w_ku-w_1v-w_2v-\cdots-w_kv$.
It follows from Theorem \ref{21t} that we obtain the desired result.
\end{pf}

\begin{cor}\label{24c}
Suppose $G_{\pi}=K_t\vee H$ with  $uv\in E(H)$,
$xy\in E(H)$,
$uy\notin E(H)$ and $xv\notin E(H)$.
Let $H'=H+uy+xv-uv-xy$ and $G'=K_t\vee H'$.
If $H'$ is connected,
then

$(i)$ $f_G(u)>f_G(x)$ if and only if $f_G(y)<f_G(v)$;
\par
$(ii)$ $f_G(u)=f_G(x)$ if and only if $f_G(y)=f_G(v)$;
Moreover,
$f_G(u)=f_G(x)$ $\big($or $f_G(y)=f_G(v)\big)$ if and only if $G'$ is also a $\Theta_\alpha$-maximal graph of $\Gamma(\pi,t;c)$.
\end{cor}

\begin{pf}
Note that $G'=G_\pi+uy+xv-uv-xy$,
and $G'\in \Gamma(\pi,t;c)$ since $H'$ is connected.
By Theorem \ref{22t} and the definition of $G_\pi$,
the result follows immediately.
\end{pf}

In order to describe the structure  of  $G_{\pi}$, namely, the  $\Theta_\alpha$-maximal graphs in the class $\Gamma(\pi,t;c)$,
we need to introduce the following concepts.

\begin{defn}\label{25d}{\em\cite{Biyikoglu,Zhang0}}
Let $G$ be a connected simple graph with    $V(G)=\big\{v_{1},v_{2},\ldots,v_n\big\}$.
We call a well-ordering $v_1\prec v_2 \prec \cdots \prec v_n$ of $V(G)$
$\big($shortly written as $\prec$$\big)$  as a breadth-first-search ordering
$\big($BFS-ordering for short$\big)$
and call $G$ a BFS-graph if $\prec$ satisfies the following two conditions:

$(i)$ $d_G(v_{1})\geq d_G(v_{2})\geq\cdots\geq d_G(v_n)$,
and $h_G(v_{1})\leq h_G(v_{2})\leq\cdots \leq h_G(v_n)$,
where $h_G(v_i)=dist_G(v_i,v_{1})$ for $i=1,2,\ldots,n$;

$(ii)$ Let $v\in N_G(u)\setminus N_G(x)$ and $y\in N_G(x)\setminus N_G(u)$ with $h_G(u)=h_G(x)=h_G(v)-1=h_G(y)-1$.
If $u\prec x$, then $v\prec y$.
 \end{defn}

\begin{defn}\label{26d}
Let $G=K_t\vee H$ be a $t$-cone graph with $n$ vertices,
where $V(H)=\big\{v_{t+1},v_{t+2},\ldots,v_n\big\}$.
If a well-ordering $v_{t+1}\prec v_{t+2}\prec\cdots\prec v_n$ of $V(H)$  satisfies the following $(i)$, $(ii)$ and $(iii)$,
then $\prec$ is called a good BFS-ordering of $V(H)$,
and $G$ is called  a good  $t$-cone BFS-graph.

$(i)$ For any two vertices $\{u,v\}\subseteq V(H)$,  if $d_G(v)>d_G(u)$, then $f_G(v)>f_G(u)$;
\par
$(ii)$ $d_G(v_{t+1})\geq d_G(v_{t+2})\geq\cdots\geq d_G(v_n)$, $f_G(v_{t+1})\geq f_G(v_{t+2})\geq\cdots\geq f_G(v_n)$,
and $h_H(v_{t+1})\leq h_H(v_{t+2})\leq\cdots \leq h_H(v_n)$,
where $h_H(v_i)=dist_H(v_i,v_{t+1})$ for $i=t+1,t+2,\ldots,n$;
\par

$(iii)$ Let $v\in N_H(u)\setminus N_H(x)$ and $y\in N_H(x)\setminus N_H(u)$ with $h_H(u)=h_H(x)=h_H(v)-1=h_H(y)-1$.
Then $f_G(u)>f_G(x)$ if and only if $f_G(y)<f_G(v)$,
and $f_G(u)=f_G(x)$ if and only if $f_G(y)=f_G(v)$.
\end{defn}

Let $G=K_t\vee H$ be a good  $t$-cone BFS-graph with $n$ vertices, and $\prec$ be a good BFS-ordering  of $V(H)$.
We use the notation $u\equiv v$ to indicate that we can interchange the positions of $u$ and $v$ in $\prec$ to obtain another good BFS-ordering of $V(H)$.

\begin{prop}\label{27p}
Let $G=K_t\vee H$ be a good $t$-cone BFS-graph. For $\big\{u_1,u_2,u_3\big\}\subseteq V(H)$,
\par
$(i)$ if $h_H(u_1)=h_H(u_2)$, then $u_1\equiv u_2$ if and only if $f_G(u_1)=f_G(u_2)$;
\par
$(ii)$ if $h_H(u_1)=h_H(u_2)=h_H(u_3)$, $u_1\equiv u_2$ and $u_2\equiv u_3$, then $u_1\equiv u_3$.
\end{prop}

\begin{pf}
We first prove $(i)$.  On one hand,
we notice that $u\prec v$ in a good BFS-ordering $\prec$ of $V(H)$ implies $f_G(u)\geq f_G(v)$.
For this reason, if $u_1\equiv u_2$,
then $u_1\prec u_2$ is in a good BFS-ordering $\prec$ of $V(H)$
and $u_2\prec' u_1$ is in another good BFS-ordering $\prec'$ of $V(H)$,
which confirms  that $f_G(u_1)\geq f_G(u_2)$ and $f_G(u_2)\geq f_G(u_1)$,
that is,
$f_G(u_1)=f_G(u_2)$.

On conversely, we are
provided  that $f_G(u_1)=f_G(u_2)$.
We interchange the positions of $u_1$ and $u_2$ in a good BFS-ordering $\prec$ of $V(H)$ to obtain a new ordering $\prec'$ of $V(H)$.
Since $f_G(u_1)=f_G(u_2)$,
we have $d_G(u_1)=d_G(u_2)$ by Definition \ref{26d} $(i)$.
Combining this with $\big\{u_1,u_2\big\}\subseteq V(H)$ and $h_H(u_1)=h_H(u_2)$,
we conclude that the ordering $\prec'$ of $V(H)$ satisfies Definition \ref{26d},
as $\prec$ satisfies Definition \ref{26d}.
Therefore,
$\prec'$ is also a good BFS-ordering of $V(H)$, and so $(i)$ holds.

Now, we turn to prove $(ii)$. In view of  $h_H(u_1)=h_H(u_2)=h_H(u_3)$,
since $u_1\equiv u_2$ and $u_2\equiv u_3$,
then by $(i)$ we get $f_G(u_1)=f_G(u_2)=f_G(u_3)$,
and again by $(i)$ we obtain the desired result.
\end{pf}

\begin{prop}\label{28p}
If $G=K_t\vee H$ is  a good $t$-cone BFS-graph, then $H$ is a BFS-graph.
\end{prop}

\begin{pf}
On one hand,
since $G=K_t\vee H$ is a good $t$-cone BFS-graph,
then there is a good BFS-ordering
$v_{t+1}\prec v_{t+2}\prec\cdots\prec v_n$ of $V(H)$.
Note that $d_H(v)=d_G(v)-t$ for each vertex $v\in V(H)$,
then by Definition \ref{26d} $(ii)$,
we have $d_H\big(v_{t+1}\big)\geq d_H\big(v_{t+2}\big)\geq\cdots\geq d_H\big(v_n\big)$,
and hence the ordering $\prec$ of
$V(H)$ satisfies Definition \ref{25d} $(i)$.

On the other hand,
to ensure the ordering $\prec$ of
$V(H)$ satisfies Definition \ref{25d} $(ii)$,
we need to prove that $v\prec y$ if $u\prec x$,
where $v\in N_H(u)\setminus N_H(x)$ and $y\in N_H(x)\setminus N_H(u)$ with $h_H(u)=h_H(x)=h_H(v)-1=h_H(y)-1$.
Since $u\prec x$ leads to $f_G(u)\geq f_G(x)$ by Definition \ref{26d} $(ii)$,
we shall consider the following two cases.

\noindent\textbf{Case 1.~} $f_G(u)>f_G(x)$.
By Definition \ref{26d} $(iii)$,
we have $f_G(v)>f_G(y)$,
which implies that $v\prec y$ by the choice of the ordering $\prec$ of $V(H)$.

\noindent\textbf{Case 2.~} $f_G(u)=f_G(x)$.
By Definition \ref{26d} $(iii)$,
we have $f_G(v)=f_G(y)$ and hence  $d_H(v)=d_H(y)$ by Definition \ref{26d} $(i)$.
Then by the fact that $h_H(v)=h_H(y)$, $f_G(v)=f_G(y)$ and Proposition \ref{27p} $(i)$,
we have  $v\equiv y$,
which means we can assert that $v\prec y$.
Otherwise,
if $y\prec v$, then we interchange the positions of $v$ and $y$ in the ordering $\prec$ of $V(H)$ to obtain a new good BFS-ordering $\prec'$ of $V(H)$.
It is easily checked that the ordering $\prec'$ of $V(H)$ satisfies Definition \ref{25d} $(i)$, and we can go forward by considering this new good BFS-ordering $\prec'$ of $V(H)$.
\end{pf}


\section{\large{Properties of $\Theta_\alpha$-maximal graphs of $\Gamma(\pi,t;c)$}}

In this section,
we study the properties of $\Theta_\alpha$-maximal graphs $G_{\pi}$ of $\Gamma(\pi,t;c)$.
If  $G=K_t\vee H\in \Gamma(\pi,t;c)$, then for any two different vertices $u$ and $v$ of $H$, since $d_H(v)=d_G(v)-t$ and $d_H(u)=d_G(u)-t$, we obtain that
\begin{equation}\label{31e}
d_G(v)>d_G(u)~~\text{iff}~~d_H(v)>d_H(u),
~~\text{and}~~d_G(v)=d_G(u)~~\text{iff}~~d_H(v)=d_H(u).
\end{equation}

\begin{lem}\label{31l}
For any two vertices $\big\{v,u\big\}\subseteq V\big(G_{\pi}\big)$,
if $d_{G_{\pi}}(v)>d_{G_{\pi}}(u)$, then $f_{G_{\pi}}(v)>f_{G_{\pi}}(u)$;
furthermore, if $f_{G_{\pi}}(v)=f_{G_{\pi}}(u)$,
then $d_{G_{\pi}}(v)=d_{G_{\pi}}(u)$.
\end{lem}

\begin{pf}
We suppose that $G_{\pi}=K_t\vee H$.
Then,
it will suffice to show that $d_{G_{\pi}}(v)>d_{G_{\pi}}(u)$ implying $f_{G_{\pi}}(v)>f_{G_{\pi}}(u)$,
as it indicates that $f_{G_{\pi}}(v)=f_{G_{\pi}}(u)$ deducing $d_{G_{\pi}}(v)=d_{G_{\pi}}(u)$.
In what follows,
by contradiction,
we assume that $u$ and $v$ are two vertices of $G_{\pi}$ with $d_{G_{\pi}}(v)>d_{G_{\pi}}(u)$,
but $f_{G_{\pi}}(v)\leq f_{G_{\pi}}(u)$.

\noindent\textbf{Case 1.~} $v\in V(K_t)$ and $u\in V(H)$.
Now by the definition of $t$-cone graph,
$uv$ is the shortest path connecting $u$ and $v$.
Denote by $W=N_{G_{\pi}}(v)\setminus \big(N_{G_{\pi}}(u)\cup \{u\}\big)$.
Since $|W|=k:=d_{G_{\pi}}(v)-d_{G_{\pi}}(u)>0$,
  we can suppose that $W=\big\{w_1,w_2,\ldots,w_k\big\}$.
Let $G_1=G_{\pi}-vw_1-vw_2-\cdots-vw_{k}+uw_1+uw_2+\cdots+uw_{k}$.
In this case, $G_1=G_{\pi}$ and hence $G_1\in \Gamma(\pi,t;c)$.
However,
Theorem \ref{21t} leads to $\Theta(G_{\pi},\alpha)<\Theta(G_1,\alpha)$,
a contradiction.

\noindent\textbf{Case 2.~} $\{v,u\}\subseteq V(H)$.
By  $(\ref{31e})$,
we have $k:=d_H(v)-d_H(u)>0$.
Since $H$ is connected,
 we can choose vertices $W=\{w_1,w_2,\ldots,w_k\}\subseteq V(H)$ such that
$W\subseteq N_H(v)\setminus N_H(u)$ and $W \cap V\big(P_H(u,v)\big)=\emptyset$.
Let $H_2=H-vw_1-vw_2-\cdots-vw_{k}+uw_1+uw_2+\cdots+uw_{k}$.
It is easily checked  that $H_2$ is connected,
and then $G_2=K_t\vee H_2\in \Gamma(\pi,t;c)$.
Then, by Corollary \ref{23c},
we get $\Theta(G_{\pi},\alpha)<\Theta(G_2,\alpha)$
contradicting the definition of $G_{\pi}$.
\end{pf}

Let $p^{(q)}$ be the $q$ copies of an integer number $p$.

\begin{lem}\label{32l}
Let $G$ be a connected simple graph with $n$ vertices,
and let $u$ and $v$ be two vertices of $G$.
If $d_G(u)=d_G(v)=n-1$,
then $f_G(u)=f_G(v)$.
\end{lem}

\begin{pf}
Note that $N_G(u)\setminus\{v\}=N_G(v)\setminus\{u\}=V(G)\setminus\{u,v\}$,
then by  (\ref{22e}) it follows that
\begin{equation}\label{32e}
\big(\Theta(G,\alpha)-\alpha(n-1)+1\big)\big( f_G(u)-f_G(v)\big)=0.
\end{equation}

 Recall that $V(G)=\big\{v_1,v_2,\ldots,v_n\big\}$.   Thus, $d(v_1)=d_1=n-1$, and  by  setting ${\bf \varphi}=\left(1,0^{(n-1)}\right)$,
it follows from the Rayleigh-Ritz theorem \cite{Horn} that
$\Theta(G,\alpha)\geq \mathcal{R}_{M_{\alpha}(G)}({\bf \varphi})=\alpha(n-1)$. Combining this with (\ref{32e}),
we have  $f_G(u)=f_G(v)$, as desired.
\end{pf}

Now, from Lemmas \ref{31l} and \ref{32l},
it immediately follows that

\begin{cor}\label{33c}
For any three vertices  $\{u,v,w\}\subseteq V\big(G_{\pi}\big)$.
If $d_{G_{\pi}}(u)=d_{G_{\pi}}(v)=n-1>d_{G_{\pi}}(w)$,
then $f_{G_{\pi}}(u)=f_{G_{\pi}}(v)>f_{G_{\pi}}(w)$.
\end{cor}

Let $G_{\pi}=K_t\vee H$.
By Corollary \ref{33c},
all vertices of $K_t$ are symmetry and have the largest $\Theta_\alpha$-weight,
and each of them is adjacent with all vertices of $H$.
So to describe the structure of $G_{\pi}$, namely,
the  $\Theta_\alpha$-maximal graphs of $\Gamma(\pi,t;c)$,
we only need to consider the structure of $H$.
In what follows,
we shall prove that
$G_\pi$ is a good $t$-cone BFS-graph.

\begin{lem}\label{34l}
Let $G_{\pi}=K_t\vee H$, where $V(H)=\big\{v_{t+1},v_{t+2},\ldots,v_n\big\}$.
Then there is a well-ordering $\prec$ of $V(H)$ satisfying Definition $\ref{26d}$ $(i)$ and $(ii)$.
\end{lem}

\begin{pf}
Clearly,
we can give an ordering $v_{t+1}\prec v_{t+2}\prec\cdots \prec v_n$ of $V(H)$
such that $v_i\prec v_j$ whenever $d_{G_{\pi}}(v_i)>d_{G_{\pi}}(v_j)$,
or $d_{G_{\pi}}(v_i)=d_{G_{\pi}}(v_j)$ and $f_{G_{\pi}}(v_i)\geq f_{G_{\pi}}(v_j)$,
where $t+1\leq i<j\leq n$.
By Lemma \ref{31l},
$d_{G_{\pi}}(v_i)>d_{G_{\pi}}(v_j)$ implies that $f_{G_{\pi}}(v_i)>f_{G_{\pi}}(v_j)$,
where $t+1\leq i<j\leq n$.
Hence $f_{G_{\pi}}(v_{t+1})\geq f_{G_{\pi}}(v_{t+2})\geq \cdots \geq f_{G_{\pi}}(v_n)$ and $d_{G_{\pi}}(v_{t+1})\geq d_{G_{\pi}}(v_{t+2})\geq\cdots\geq d_{G_{\pi}}(v_n)$.

To complete the proof of this result, it suffices to  prove that $h_H(v_{i})\leq h_H(v_{i+1})$ holds for any $i=t+1,t+2,\ldots,n-1$
by induction on $i$.
Obviously,
for $i=t+1$,
we have $h_H(v_{t+1})=0\leq h_H(v_{t+2})$ and hence  the assertion holds.
Now we may assume that the assertion already holds for $t+1\leq i\leq k-1$,
namely,
we already have $h_H(v_{t+1})\leq\cdots \leq h_H(v_{k-1})\leq h_H(v_{k})$,
and we  will  prove $h_H(v_{k})\leq h_H(v_{k+1})$,
where $t+2\leq k\leq n-1$.
We consider the following two cases.

\noindent\textbf{Case 1.~} $f_{G_{\pi}}(v_k)>f_{G_{\pi}}(v_{k+1})$.

Since $H$ is connected,
there exists the smallest integer $s\in \big\{t+1,t+2,\ldots,k\big\}$
such that $v_s$ is adjacent to some vertex \big(say $v_q$\big) in $\big\{v_{k+1},v_{k+2},\ldots,v_n\big\}$.

\par\medskip

\noindent\textbf{Claim 1.~} $h_H(v_s)\geq h_H(v_k)-1$.

By contradiction,
suppose that $h_H(v_s)<h_H(v_k)-1$.
Thus $s\neq k$,
and it leads to $t+1\leq s\leq k-1$.
Moreover,
$v_sv_k\notin E(H)$
\big(otherwise, $h_H(v_s)\geq h_H(v_k)-1$,
contradicting our hypothesis\big).
Let $p=\min\big\{j~|~v_j\in N_H(v_k)\big\}$.
If $v_p\prec v_s$ or $v_p=v_s$,
then $t+1\leq p\leq k-1$,
and bearing in mind that $v_p\in N_H(v_k)$,
it follows from the induction hypothesis that
$h_H(v_k)\leq h_H(v_p)+1\leq h_H(v_s)+1<h_H(v_k)$,
a contradiction.
Hence $v_s\prec v_p$,
and we get $f_{G_\pi}(v_s)\geq f_{G_\pi}(v_p)\geq f_{G_\pi}(z)$ for every $z\in N_H(v_k)$ by the choice of $p$.
Since $s\leq k-1$ and $q\geq k+1$,
we have $f_{G_{\pi}}(v_s)\geq f_{G_{\pi}}(v_k)>f_{G_{\pi}}(v_{k+1})\geq f_{G_{\pi}}(v_q)$.
Note that $H$ is connected and $v_sv_k\notin E(H)$,
there exists a shortest path,
say $P_H(v_s,v_k)=v_s\cdots yv_k$ from $v_s$ to $v_k$ in $H$,
where $v_s\neq y\neq v_k$.
Besides, since $f_{G_{\pi}}(v_k)>f_{G_{\pi}}(v_q)$,
by Lemma \ref{31l} we have $d_{G_{\pi}}(v_k)\geq d_{G_{\pi}}(v_q)$,
and then $d_H(v_k)\geq d_H(v_q)$ by $(\ref{31e})$.

In what follows,
we will prove that there exists $z\in N_H(v_k)\setminus N_H(v_q)$ such that
$H'=H+v_sv_k+v_qz-v_sv_q-v_kz$ is connected,
and then $G':=K_t\vee H'=G_{\pi}+v_sv_k+v_qz-v_sv_q-v_kz \in \Gamma(\pi,t;c)$.
If $v_sv_q\in E\big(P_H(v_s,v_k)\big)$, then
since $d_H(v_k)\geq d_H(v_q)$ and $v_s\in N_H(v_q)\setminus N_H(v_k)$,
there must exist vertex $z\in N_H(v_k)\setminus V\big(P_H(v_s,v_k)\big)$ such that
$v_qz\notin E(H)$, as desired.
Now we suppose that   $v_sv_q\notin E\big(P_H(v_s,v_k)\big)$. If $v_qy\notin E(H)$, then
  $z=y$ as desired. Otherwise, $v_qy\in E(H)$.
In such case,
again by the fact that $d_H(v_k)\geq d_H(v_q)$ and $v_s\in N_H(v_q)\setminus N_H(v_k)$,
there must exist vertex $z\in N_H(v_k)\setminus V(P_H(v_s,v_k))$ such that
$v_qz\notin E(H)$, as desired.

Since we have proved that $f_{G_{\pi}}(v_s)\geq f_{G_{\pi}}(v_p)\geq f_{G_{\pi}}(z)$
for any $z\in N_H(v_k)$ and $f_{G_{\pi}}(v_k)>f_{G_{\pi}}(v_q)$,
it follows from Theorem \ref{22t} that $\Theta(G',\alpha)>\Theta(G_{\pi},\alpha)$,
contradicting the definition  of  $G_{\pi}$.
Therefore, Claim 1 holds.

\par\medskip

Let $P_H\big(v_{t+1},v_{k+1}\big)$ be a shortest path from $v_{t+1}$ to $v_{k+1}$ in $H$, and let $v_r$ be the last vertex belonging to $\big\{v_{t+1},v_{t+2},\ldots,v_{k}\big\}$ on the path $P_H\big(v_{t+1},v_{k+1}\big)$.
If $s\leq k-1$,
then $h_H(v_s)\leq h_H(v_r)$ by the induction hypothesis and the choice of $v_s$.
If $s=k$,
by the choice of $v_s$ we have $r=s=k$,
and then $h_H(v_s)=h_H(v_r)$.
Combining these with Claim 1,
we obtain that
$$h_H(v_{k+1})\geq h_H(v_r)+1\geq h_H(v_s)+1\geq h_H(v_k).$$

\noindent\textbf{Case 2.~} $f_{G_{\pi}}(v_k)=f_{G_{\pi}}(v_{k+1})$.

Let $d=\min\big\{i~|~f_{G_{\pi}}(v_i)=f_{G_{\pi}}(v_{k+1})~\text{and}~i\geq t+1\big\}$.
Then $t+1\leq d\leq k$ and $f_{G_{\pi}}(v_d)=f_{G_{\pi}}(v_{d+1})=\cdots=f_{G_{\pi}}(v_k)=f_{G_{\pi}}(v_{k+1})$.
By Lemma \ref{31l},
we have $d_{G_{\pi}}(v_d)=d_{G_{\pi}}(v_{d+1})=\cdots=d_{G_{\pi}}(v_k)=d_{G_{\pi}}(v_{k+1})$.
At this time,
we relabel $V(H)$ by the following method:
If $j\notin\big\{d,d+1,\ldots,k,k+1\big\}$, then let $v_j'=v_j$;
if $j\in \{d,d+1,\ldots,k,k+1\big\}$, then let $v_{d+p-1}'=v_j$ if and only if
$h_H(v_j)$ is the $p$-th smallest value among $\big\{h_H(v_d),h_H(v_{d+1}),\ldots,h_H(v_k),h_H(v_{k+1})\big\}$.
Then by the induction hypothesis,
there is an ordering $v'_{t+1}\prec v'_{t+2}\prec\cdots \prec v'_n$ of
$V(H)=\big\{v'_{t+1},v'_{t+2},\ldots,v'_n\big\}$ such that
$d_H\big(v'_{t+1}\big)\geq d_H\big(v'_{t+2}\big)\geq\cdots\geq d_H\big(v'_n\big),$
$f_{G_{\pi}}\big(v_{t+1}'\big)\geq f_{G_{\pi}}\big(v_{t+2}'\big)\geq \cdots \geq f_{G_{\pi}}\big(v_n'\big),$
and
$$h_H\big(v'_{t+1}\big)\leq h_H\big(v'_{t+2}\big)\leq \cdots\leq h_H\big(v'_{d-1}),~~h_H(v'_{d}\big)\leq h_H\big(v'_{d+1}\big)\leq \cdots\leq h_H\big(v'_{k}\big)\leq h_H\big(v'_{k+1}\big).$$
Clearly,
if $d=t+1$,
the result already holds;
if $d\geq t+2$,
we get $f_{G_{\pi}}\big(v'_{d-1}\big)>f_{G_{\pi}}\big(v'_d\big)$ by the choice of $d$,
and analogously as the proof of Case 1 we obtain that
$h_H\big(v'_{d-1}\big)\leq h_H\big(v'_{d}\big)$.
Therefore,
$$h_H\big(v'_{t+1}\big)\leq h_H\big(v'_{t+2}\big)\leq \cdots\leq h_H\big(v'_{d-1}\big)\leq h_H\big(v'_{d}\big)\leq h_H\big(v'_{d+1}\big)\leq \cdots\leq h_H\big(v'_{k}\big)\leq h_H\big(v'_{k+1}\big).$$
This completes the proof.
\end{pf}

\begin{thm}\label{35t}
If $G_{\pi}=K_t\vee H$,
then $G_\pi$ is a good $t$-cone BFS-graph
and $H$ is a BFS-graph.
\end{thm}

\begin{pf}
By Lemma \ref{34l},
there exists a well-ordering   $\prec$  of
$V(H)$
satisfying Definition \ref{26d} $(i)$ and $(ii)$.
Now it suffices to show that the ordering $\prec$ of
$V(H)$ satisfies Definition \ref{26d} $(iii)$.
By the condition that $v\in N_H(u)\setminus N_H(x)$ and $y\in N_H(x)\setminus N_H(u)$ with $h_H(u)=h_H(x)=h_H(v)-1=h_H(y)-1$ in Definition \ref{26d} $(iii)$,
let $H'=H+uy+xv-uv-xy$ and $G'=K_t \vee H'$.
It can be checked that $H'$ is connected.
Thus,
Corollary \ref{24c} implies that Definition \ref{26d} $(iii)$ holds.
Hence $G_\pi$ is a good $t$-cone BFS-graph,
and so $H$ is a BFS-graph by Proposition \ref{28p}.
\end{pf}

\begin{lem}\label{36l}
Let $G_{\pi}=K_t\vee H$.
Let $uv$ be an edge on a cycle of $H$
and $P=wx_1x_2\cdots x_{s}$ $(s\geq 1)$ be a path of $H$ such that $f_{G_{\pi}}(x_s)<\min\big\{f_{G_{\pi}}(u),f_{G_{\pi}}(v)\big\}$.
Suppose that $N_H(v)\cap \big\{w,x_2,x_4,x_6,\ldots,x_{q_1}\big\}=\emptyset$
and $N_H(u)\cap \big\{x_1,x_3,x_5,\ldots,x_{q_2}\big\}=\emptyset$,
 where $q_1+1=s=q_2$ if $s$ is odd
and $q_1=s=q_2+1$ if $s$ is even.
Then, $f_{G_{\pi}}(u)>f_{G_{\pi}}(w)$.
\end{lem}

\begin{pf}
By contradiction, suppose that $f_{G_{\pi}}(u)\leq f_{G_{\pi}}(w)$.
Let $H'=H+vw+ux_1-uv-wx_1$ and $G'=K_t\vee H'$.
Since $uv$ is an edge on a cycle of $H$,
then $H'$ is connected.
In view of  the definition  of  $G_{\pi}$,
by Corollary \ref{24c} and the assumption $f_{G_{\pi}}(u)\leq f_{G_{\pi}}(w)$,
we get $f_{G_{\pi}}(v)\leq f_{G_{\pi}}(x_1)$.
Let $H''=H+ux_1+vx_2-uv-x_1x_2$ and $G''=K_t\vee H''$.
Clearly,
$H''$ is connected by the fact that $uv$ being an edge on a cycle of $H$.
Since $f_{G_{\pi}}(v)\leq f_{G_{\pi}}(x_1)$.
It follows from Corollary \ref{24c} and the definition  of $G_{\pi}$ that $f_{G_{\pi}}(u)\leq f_{G_{\pi}}(x_2)$.
By repeating similar arguments,
we finally deduce  that $\min\big\{f_{G_{\pi}}(u),f_{G_{\pi}}(v)\big\}\leq f_{G_{\pi}}(x_s)$,
contradicting  the condition $f_{G_{\pi}}(x_s)<\min\big\{f_{G_{\pi}}(u),f_{G_{\pi}}(v)\big\}$.
\end{pf}

Let $G=K_t\vee H$ be a $t$-cone $c$-cyclic graph with $n$ vertices.
If $c\geq 1$,
then we define the {\em basic graph} of $H$,
denoted by $\mathcal{B}(H)$,
as the graph obtained from $H$ by recursively deleting pendant vertices
(that is, those vertices of degree one),
to the resultant graph until no pendant vertices remain.
From the definition of basic graph,
$\mathcal{B}(H)$ is unique and it is also a $c$-cyclic graph.

\begin{thm}\label{37t}
Let $G_{\pi}=K_t\vee H$.
If  $u\in V\big(\mathcal{B}(H)\big)$ and $w\in V(H)\setminus V\big(\mathcal{B}(H)\big)$,
then $f_{G_{\pi}}(u)>f_{G_{\pi}}(w)$.
\end{thm}

\begin{pf}
 Note that $u\in V\big(\mathcal{B}(H)\big)$ and $w\in V(H)\setminus V\big(\mathcal{B}(H)\big)$.
If $d_H(w)=1$,
since $d_H(u)\geq 2$,
then $f_{G_{\pi}}(u)>f_{G_{\pi}}(w)$ by $(\ref{31e})$ and Lemma \ref{31l}.
In the following,
we assume that $d_H(w)\geq 2$.

\noindent\textbf{Case 1.~} $u$ lies on some cycle $C$ of $H$.
By the definition of $\mathcal{B}(H)$ and since $w\not\in V\big(\mathcal{B}(H)\big)$,
there exists some pendant path in $H$,
say $P=wx_1x_2\cdots x_s$ $(s\geq 1)$,
where $d_H(x_s)=1$,
$V(P)\cap V\big(\mathcal{B}(H)\big)=\emptyset$,
and $N_H(u)\cap \big\{x_1,x_2,\ldots,x_s\big\}=\emptyset$.
By the choice of $P$,
there exists an edge $uv$ in  the cycle $C$
such that $vw\notin E(H)$ and
$\big(N_H(u)\cup N_H(v)\big)\cap \big\{x_1,x_2,\ldots,x_s\big\}=\emptyset$.
Moreover,
since $\min\big\{d_H(u),d_H(v)\big\}\geq 2$ and $d_H(x_s)=1$,
by $(\ref{31e})$ and Lemma \ref{31l}
we get $\min\big\{f_{G_{\pi}}(u),f_{G_{\pi}}(v)\big\}>f_{G_{\pi}}(x_s)$, and so  $f_{G_{\pi}}(u)>f_{G_{\pi}}(w)$ follows from  Lemma \ref{36l}.

\noindent\textbf{Case 2.~} $u$ does not lie on any cycle of $H$.
In such case,
$u$ lies on a path $P'$ of $H$,
where $P'$ is the unique path connecting two cycles, say $C_1$ and $C_2$, of $H$.
Suppose that $\{x\}=V(P')\cap V(C_1)$ and $\{y\}=V(P')\cap V(C_2)$.
Notice that there is a good BFS-ordering $v_{t+1}\prec v_{t+2}\prec \cdots \prec v_{n}$ of $V(H)$ by Theorem \ref{35t},
we consider the following two subcases.

\noindent\textbf{Subcase 2.1~} $u\in V\big(P_H(v_{t+1},x)\big)$ or $u\in V\big(P_H(v_{t+1},y)\big)$. By the symmetry of $x$ and $y$, without loss of generality, we suppose that $u\in V\big(P_H(v_{t+1},x)\big)$.
Then by Theorem \ref{35t},
we have $f_{G_{\pi}}(u)\geq f_{G_{\pi}}(x)$.
Moreover,
since $x$ lies on cycle $C_1$,
similarly as the proof of Case 1,
we have $f_{G_{\pi}}(x)>f_{G_{\pi}}(w)$ and so $f_{G_{\pi}}(u)\geq f_{G_{\pi}}(x)>f_{G_{\pi}}(w)$, as desired.

\noindent\textbf{Subcase 2.2~} $u\notin V\big(P_H(v_{t+1},x)\big)\cup  V\big(P_H(v_{t+1},y)\big)$.
In such subcase,
it is easily checked  that $u$ lies on some cycle of $H$,
which is a contradiction.
\end{pf}

\begin{thm}\label{38t}
Let $G_{\pi}=K_t\vee H$   with  $\big\{u,v\big\}\subseteq V\big(\mathcal{B}(H)\big)$.
If $d_{\mathcal{B}(H)}(u)>d_{\mathcal{B}(H)}(v)$,
then $f_{G_{\pi}}(u)>f_{G_{\pi}}(v)$.
Moreover,
$f_{G_{\pi}}(u)=f_{G_{\pi}}(v)$ implies that $d_{\mathcal{B}(H)}(u)=d_{\mathcal{B}(H)}(v)$
\end{thm}

\begin{pf}
It suffices to show that
$d_{\mathcal{B}(H)}(u)>d_{\mathcal{B}(H)}(v)$ implying $f_{G_{\pi}}(u)>f_{G_{\pi}}(v)$,
as this  indicates that $f_{G_{\pi}}(u)=f_{G_{\pi}}(v)$ deducing $d_{\mathcal{B}(H)}(u)=d_{\mathcal{B}(H)}(v)$.
If $d_H(u)>d_H(v)$,
then by $(\ref{31e})$ and Lemma \ref{31l},
we obtain that $f_{G_{\pi}}(u)>f_{G_{\pi}}(v)$.
So we may assume that $d_H(u)\leq d_H(v)$ in the following.
Combining this with $d_{\mathcal{B}(H)}(u)>d_{\mathcal{B}(H)}(v)$,
we conclude that there exists some vertex $z$ such that $z\in N_H(v)\setminus N_H(u)$ and $z\in V(H)\setminus V(\mathcal{B}(H))$.
On the other hand,
again by the fact that $d_{\mathcal{B}(H)}(u)>d_{\mathcal{B}(H)}(v)$,
there exists some vertex $w\in N_{\mathcal{B}(H)}(u)\setminus N_{\mathcal{B}(H)}(v)$ such that $w\notin V\big(P_{\mathcal{B}(H)}(u,v)\big)$.
By Theorem \ref{37t},
we have $f_{G_{\pi}}(w)>f_{G_{\pi}}(z)$.
Let $H'=H+uz+vw-uw-vz$ and $G'=K_t\vee H'$.
Clearly,
$H'$ is connected.
Hence by Corollary \ref{24c},
we get $f_{G_{\pi}}(u)>f_{G_{\pi}}(v)$.
\end{pf}

An \emph{internal path} joining $u_1$ and $u_{k+1}$
\big(need not be distinct\big) of $G$
is a path $Q=u_1u_2\cdots u_{k+1}$ such that
$d_G(u_1)\geq 3$, $d_G(u_{k+1})\geq 3$,
and $d_G(u_i)=2$ \big(if exist\big) for $2\leq i\leq k$.
Let $T$ be a tree such that $V(G)\cap V(T)=\emptyset$.
If we obtain a new graph $G'$ from $T$ and $G$ by identifying
one vertex $v$ of $T$ and some vertex of $G$,
then $T$ is called a \emph{root tree} of $G'$,
and $v$ is called the \emph{root} of $T$.

\begin{lem}\label{39l}
Let ${G_{\pi}}=K_t\vee H$ and $Q=u_1u_2\cdots u_{k+1}$ $\big(k\geq 1\big)$ be an internal path of $H$ $\big($or $\mathcal{B}(H)$$\big)$.
If $k\geq 2$ and $u_1\neq u_{k+1}$, then $Q$ is a part of some cycle of $H$. Furthermore, if $k=2$, $u_1\neq u_{3}$ and  there is a pendant vertex not on the root tree  with root $u_2$, then $u_1u_{3}\in E(H).$
\end{lem}

\begin{pf}
Note that $V(H)$ has a good BFS-ordering $v_{t+1}\prec v_{t+2}\prec \cdots \prec v_{n}$ by Theorem \ref{35t}.
By contradiction,
we assume that $Q$ is not  a part of any  cycle of $H$
when $k\geq 2$ and $u_1\neq u_{k+1}$.
Then, $Q$ is the unique path of $H$ connecting $u_1$ and
$u_{k+1}$.
 Combining this with the connectivity of $H$,
either $u_2\in P_{H}\big(v_{t+1},u_1\big)$
or  $u_2\in P_{H}\big(v_{t+1},u_{k+1}\big)$.
Then by Theorem \ref{35t},
$f_{G_{\pi}}(u_2)\geq \min\big\{f_{G_{\pi}}(u_1),f_{G_{\pi}}(u_{k+1})\big\}$.  On the other hand,
since $Q$ is an internal path of $H$ \big(or $\mathcal{B}(H)$\big),
we have $d_H(u_2)=2<3\leq \min\big\{d_{H}(u_1),d_{H}(u_{k+1})\big\}$
\big(or $d_{\mathcal{B}(H)}(u_2)=2<3\leq \min\big\{d_{\mathcal{B}(H)}(u_1),d_{\mathcal{B}(H)}(u_{k+1})\big\}$\big),
and it follows from $(\ref{31e})$ and Lemma \ref{31l} \big(or Theorem \ref{38t}\big) that
$f_{G_{\pi}}(u_2)<\min\big\{f_{G_{\pi}}(u_1),f_{G_{\pi}}(u_{k+1})\big\}$,
a contradiction.
Now, we can conclude that $Q$ is a part of some cycle of $H$.

Next, we turn to prove the `Furthermore' part.
Now $k=2$, $u_1\neq u_3$,
and we suppose that $u_1u_3\not\in E(H)$ by contradiction.
Let $x$ be a pendant vertex of $H$ not on the root tree with root $u_2$. Without loss of generality, we may suppose that $dist_H(u_1,x)\leq dist_H(u_3,x)$ by symmetry.
Since $d_H(x)=1$ and $u_1\neq u_3$,
we can suppose that $P_H(u_1,x)=wx_1x_2\cdots x_{s}$
by setting $w=u_1$ and $x_s=x$,
where $s\geq 1$.
Moreover,
bearing in mind that $dist_H(u_1,x)\leq dist_H(u_3,x)$ and $u_1u_3\notin E(H)$,
we have $N_H(u_3)\cap \big(P_H(u_1,x)\setminus\{x_1\}\big)=\emptyset$,
that is, $N_H(u_3)\cap \big\{w,x_2,x_3,\ldots,x_s\big\}=\emptyset$.
According to the definitions of $Q$ and $x$,
we have $N_H(u_2)\cap \big\{x_1,x_2,\ldots,x_s\big\}=\emptyset$.
Note that $\min\big\{d_{H}(u_2),d_{H}(u_{3})\big\}\geq 2>d_{H}(x)=1$,
thus $\min\big\{f_{G_{\pi}}(u_2),f_{G_{\pi}}(u_{3})\big\}>f_{G_{\pi}}(x)$
by $(\ref{31e})$ and Lemma \ref{31l}.
Then by taking $uv=u_2u_3$ with $v=u_3$,
and $P=P_H(u_1,x)=wx_1x_2\cdots x_{s}$ in Lemma \ref{36l},
we have $f_{G_{\pi}}(u_2)>f_{G_{\pi}}(u_1)$.
However,
since $d_H(u_1)\geq 3>2=d_H(u_2)$
\big(or $d_{\mathcal{B}(H)}(u_1)\geq 3>2=d_{\mathcal{B}(H)}(u_2)$\big)
by the fact that $Q$ is an internal path of $H$ \big(or $\mathcal{B}(H)$\big),
then by $(\ref{31e})$ and Lemma \ref{31l} \big(or Theorem \ref{38t}\big),
we get $f_{G_{\pi}}(u_2)<f_{G_{\pi}}(u_1)$,
a contradiction.
\end{pf}

\begin{thm}\label{310t}
Let ${G_{\pi}}=K_t\vee H$ with $\delta(H)=1$,
and $Q=u_1u_2\cdots u_{k+1}$ $\big(k\geq 1\big)$ be an internal path of $H$.
$(i)$ If $u_1=u_{k+1}$, then $k=3$;
$(ii)$ If $u_1\neq u_{k+1}$,
then $k\leq 2$,
and $u_1u_3\in E(H)$ when $k=2$.
\end{thm}

\begin{pf}
We shall prove the results by contradiction,
that is,
we may assume that $k\geq 4$ when $u_1=u_{k+1}$,
and $k\geq 3$ when $u_1\neq u_{k+1}$.
 No matter which case happens,
we always have
$d_{H}(u_{2})=d_{H}(u_{3})=\cdots=d_{H}(u_{k})=2$,
and $Q$ is a part of some cycle of $H$
\big(this conclusion is inevitable in the former case,
and can be deduced  from Lemma \ref{39l} in the latter case\big).
 So, $V(\mathcal{B}(H))\neq \emptyset$.
Since $\delta(H)=1$,
there exists some pendant path in $H$,
say $P=wx_1x_2\cdots x_s$ $\big(s\geq 1\big)$,
such that $d_H(x_s)=1$ and $V(P)\cap V\big(\mathcal{B}(H)\big)=\big\{w\big\}$.
Clearly, $d_H(w)\geq 3$.
Moreover,
by the definition of $Q$,
either $u_2w\not\in E(H)$ or $u_3w\not\in E(H)$,
 $w\not\in \big\{u_2,u_3\big\}$, and
$\big(N_H(u_2)\cup N_H(u_3)\big)\cap \big\{x_1,x_2,\ldots,x_s\big\}=\emptyset$.
Without loss of generality,
we may suppose that $u_3w\not\in E(H)$ by the symmetry of $u_2$ and $u_3$.
Note that $\min\big\{d_{H}(u_2),d_{H}(u_{3})\big\}=2>d_{H}(x_s)=1$,
thus $\min\big\{f_{G_{\pi}}(u_2),f_{G_{\pi}}(u_{3})\big\}>f_{G_{\pi}}(x_s)$
by $(\ref{31e})$ and Lemma \ref{31l}.
By considering  the edge $uv=u_2u_3$ with $v=u_3$ in Lemma \ref{36l}, we will deduce that $f_{G_{\pi}}(u_2)>f_{G_{\pi}}(w)$,
contradicting to Lemma \ref{31l} and $(\ref{31e})$, as $d_H(w)\geq 3>d_H(u_2)=2$. Thus, $k\leq 3$ when $u_1=u_{k+1}$
\big(there must be $k=3$ since $G_\pi$ is simple\big),
and $k\leq 2$ when $u_1\neq u_{k+1}$.

If $k=2$ and $u_1\neq u_3$,
then since $d_H(u_2)=2$ and $\delta(H)=1$,  there is a pendant vertex not on the root tree  with root $u_2$, and so  $u_1u_{3}\in E(H)$ by Lemma \ref{39l}.
\end{pf}

\begin{thm}\label{311t}
Let $G_{\pi}=K_t\vee H$ with $\delta(H)=1$, and
$Q=u_1u_2\cdots u_{k+1}$ $\big(k\geq 1\big)$ be an internal path of $\mathcal{B}(H)$.
$(i)$ If $u_1=u_{k+1}$, then $k=3$;
$(ii)$ If $u_1\neq u_{k+1}$, then $k\leq 2$;
furthermore,
when $k=2$,
then either $u_1u_3\in E(H)$ or all the pendant vertices of $H$
are on the root tree with root $u_2$.
\end{thm}

\begin{pf}
 Clearly, $k\geq 3$ when $u_1=u_{k+1}$.
By contradiction,
we will prove the results by assuming that $k\geq 4$ when $u_1=u_{k+1}$,
and $k\geq 3$ when $u_1\neq u_{k+1}$.
No matter which case happens,
$Q$ is a part of some cycle of $H$
(this conclusion is inevitable in the former case,
and can be derive from Lemma \ref{39l} in the latter case).

If $d_H(u_2)=2$,
then we suppose that $d_H(u_2)=d_H(u_3)=\cdots=d_H(u_q)=2<3\leq d_H(u_{q+1})$,   where $2\leq q\leq k$.
Then,
$Q'=u_1u_{2}u_{3}\cdots u_{q+1}$ is an internal path of $H$.
Since we assume that $k\geq 4$ when $u_1=u_{k+1}$
and $k\geq 3$ when $u_1\neq u_{k+1}$,
we can deduce that $2\leq q <k$ by Theorem \ref{310t}.
Once again,
since $u_1\neq u_{q+1}$ by $2\leq q <k$
and $Q'$ is an internal path of $H$,
 by Theorem \ref{310t} $(ii)$ we deduce that
$q=2$ and $u_1u_3\in E(H)$.
 Note that $Q$ is an internal path of $\mathcal{B}(H)$.
On one hand, we have $d_{\mathcal{B}(H)}\big(u_{j}\big)=2$ for $j\in \big\{2,3,\ldots,k\big\}$.
On the other hand, $\big\{u_1,u_2,u_4\big\}\subseteq N_{\mathcal{B}(H)}(u_3)$
and then $d_{\mathcal{B}(H)}\big(u_{3}\big)\geq 3$,
a contradiction.
Thus, $d_H(u_2)\geq 3$.

We also claim that $d_H(u_3)\geq 3$.
Otherwise,
we have $d_H(u_3)=2$,
then we suppose that $d_H(u_3)=d_H(u_4)=\cdots=d_H(u_q)=2<3\leq d_H(u_{q+1})$,  where $3\leq q\leq k$,
and this implies that $Q''=u_{2}u_{3}\cdots u_{q+1}$ is an internal path of $H$.
Since $u_{2}\neq u_{q+1}$,
  it follows from   Theorem \ref{310t} $(ii)$  that  $q=3$ and $u_2u_{4}\in E(H)$.
Thus,
$\big\{u_1,u_3,u_4\big\} \subseteq N_{\mathcal{B}(H)}(u_2)$
and so $d_{\mathcal{B}(H)}\big(u_{2}\big)\geq 3$,
a contradiction.
This contradiction confirms our claim.
With the similar reason,
we can conclude that $d_H(u_j)\geq 3$,
that is,
$u_j$ is a root of some root tree for any $j\in \big\{2,3,\ldots,k\big\}$.

Now, we let $x$ \big(resp., $y$\big)  be a vertex in the neighbor set of $u_2$ \big(resp., $u_3$\big)  such that $x$ \big(resp., $y$\big) is in the root tree   of $H$ with root $u_2$ \big(resp., $u_3$\big).
It is easily checked that $\big\{x\big\}\cap \big(N_H(u_3)\cup N_H(u_4)\big)=\emptyset$ and
$\big\{y\big\}\cap \big(N_H(u_1)\cup N_H(u_2)\big)=\emptyset$.
On one hand, since $\big\{u_3,u_4\big\}\subseteq \mathcal{B}(H)$ and $x\in V(H)\setminus V\big(\mathcal{B}(H)\big)$, we have  $\min\big\{f_{G_{\pi}}(u_3),f_{G_{\pi}}(u_4)\big\}>f_{G_{\pi}}(x)$ by Theorem \ref{37t},
and  hence $f_{G_{\pi}}(u_3)>f_{G_{\pi}}(u_2)$ by Lemma \ref{36l}, as $u_4 u_2\not\in E(H)$.
On the other hand, it follows from Theorem \ref{37t} that  $\min\big\{f_{G_{\pi}}(u_1),f_{G_{\pi}}(u_2)\big\}>f_{G_{\pi}}(y)$ \big(since $\big\{u_1,u_2\big\}\subseteq \mathcal{B}(H)$ and $y\in V(H)\setminus V\big(\mathcal{B}(H)\big)$\big),
then by Lemma \ref{36l} we have $f_{G_{\pi}}(u_2)>f_{G_{\pi}}(u_3)$, as  $u_1 u_3\not\in E(H)$.
It is a contradiction.

All in all,
we deduce that $k=3$ when $u_1=u_{k+1}$,
and $k\leq 2$ when $u_1\neq u_{k+1}$.
Finally,
if $k=2$,
$u_1\neq u_3$ and there exists at least one pendant vertex  not on the root tree  with root $u_2$,
then by Lemma \ref{39l},
we have $u_1u_3\in E(H)$,
and this completes the proof.
\end{pf}


\section{\large{The unique  $\Theta_\alpha$-maximal graph  of $\Gamma(\pi,t;c)$ for $c\in \big\{0,1,2\big\}$}}

In this section, we will prove that the $\Theta_\alpha$-maximal graph  of $\Gamma(\pi,t;c)$ is unique for $c\in \big\{0,1,2\big\}$.
Recall that $\pi=\big(d_1,d_2,\ldots,d_n\big)$ is a non-increasing degree sequence of $t$-cone $c$-cyclic graph,
and $\pi^*=\big(d^*_{t+1},d^*_{t+2},\ldots,d^*_{n}\big)$,
where $d^*_j=d_j-t$ for $j\in \big\{t+1,t+2,\ldots,n\big\}$.
Firstly, we have the following assertion on $\pi$ and $\pi^*$, which follows from the   definition of $t$-cone $c$-cyclic graph.

\begin{prop}\label{41p}
If $0\leq c\leq {n-t \choose 2}-n+t+1$, then
$d_1=d_2=\cdots=d_t=n-1$ and $\sum_{i=t+1}^n  d^{*}_i =2(n-t+c-1).$ Moreover, we have \par

$(i)$ if $c=0$,
then

\noindent $(4.0.1)~~0\leq t\leq n-2~~\text{and}~~d^*_{t+1}\geq d^*_{t+2}\geq \cdots\geq d^*_n=1;$

\par $(ii)$ if $c=1$,
then one of the following holds:

\noindent$(4.1.1)~~0\leq t\leq n-3~~\text{and}~~d^*_{t+1}=d^*_{t+2}= \cdots= d^*_n=2;$

\noindent$(4.1.2)~~0\leq t\leq n-4,~~d^*_{t+1}\geq 3,~~d^*_{t+2}\geq d^*_{t+3}\geq 2~~\text{and}~~d^*_{t+4}\geq d^*_{t+5}\geq \cdots\geq d^*_n=1;$ \par
$(iii)$ if $c=2$,
then one of the following holds:

\noindent$(4.2.1)~~0\leq t\leq n-5,~~d^*_{t+1}=4~~\text{and}~~d^*_{t+2}=d^*_{t+3}= \cdots= d^*_n=2;$

\noindent$(4.2.2)~~0\leq t\leq n-4,~~d^*_{t+1}=d^*_{t+2}=3~~\text{and}~~d^*_{t+3}=d^*_{t+4}= \cdots= d^*_n=2;$

\noindent$(4.2.3)~~0\leq t\leq n-6,~~d^*_{t+1}\geq 5,~~d^*_{t+2}=d^*_{t+3}=d^*_{t+4}=d^*_{t+5}=2\geq d^*_{t+6}\geq d^*_{t+7}\geq \cdots\geq d^*_n=1;$

\noindent$(4.2.4)~~0\leq t\leq n-5,~~d^*_{t+1}\geq d^*_{t+2}\geq 3,~~d^*_{t+3}\geq d^*_{t+4}\geq 2,~~\text{and}~~d^*_{t+5}\geq d^*_{t+6}\geq \cdots\geq d^*_n=1.$
\end{prop}

Before starting our main results,
we need to introduce a special $t$-cone tree
(resp., $t$-cone unicyclic graph and $t$-cone bicyclic graph) as follows.

When
$\pi=\big(d_1,d_2,\ldots,d_n\big)$ is the degree sequence of a $t$-cone tree, then $\pi^{*}=\big(d^*_{t+1},d^*_{t+2},\ldots,d^*_{n}\big)$ is a tree degree sequence, and we use the following breadth-first-search method to construct a tree  $T_{\pi^{*}}$  with $\pi^{*}$ as its degree sequence:

$(i)$ Put $s_0=1$, select a vertex $v_{0,1}$ as a root of $T_{\pi^{*}}$ and begin with $v_{0,1}$ of the zeroth layer;

$(ii)$ Put $s_1=d^*_{t+1}$, select $s_1$ vertices $\big\{v_{1,1},v_{1,2},\ldots,v_{1,s_1}\big\}$ of the first layer such that they are adjacent to $v_{0,1}$, and so   $d_{T_{\pi^{*}}}(v_{0,1})=s_1=d^*_{t+1}$;

$(iii)$  The remaining vertices of $T_{\pi^*}$ appear in a {\em BFS-connecting}, that is:
assume that all vertices of the $p$-th layer of $T_{\pi^{*}}$ have been constructed and they are $\big\{v_{p,1},v_{p,2},\ldots,v_{p,s_p}\big\}$;
by using the induction hypothesis,
now we construct all the vertices of the $(p+1)$-st layer;
let $d_{T_{\pi^{*}}}(v_{p,i})=d^*_{t+i+\sum_{j=0}^{p-1}s_j}$ for $i=1,2,\ldots,s_p$
\big(here, we shall stop if $t+i+\sum_{j=0}^{p-1}s_j>n$\big);
select $s_{p+1}=\sum_{i=1}^{s_p} d_{T_{\pi^{*}}}(v_{p,i})-s_p$ vertices $\big\{v_{p+1,1},v_{p+1,2},\ldots,v_{p+1,s_{p+1}}\big\}$ for the $(p+1)$-st layer such that
$d_{T_{\pi^{*}}}(v_{p,i})-1$ vertices are adjacent to $v_{p,i}$ for $i=1,2,\ldots,s_p$.

In this way, we obtain a unique  tree $T_{\pi^{*}}$ with $\pi^*$ as its degree sequence, and it can be easily checked   that $T_{\pi^{*}}$ is a BFS-tree.
Furthermore, for any given tree degree sequence $\pi^*$,
Zhang \cite{Zhang0} had  shown    that the
corresponding  BFS-tree is unique  and it is isomorphic to $T_{\pi^{*}}$.
Let $T_{\pi}^{(t)}=K_t \vee T_{\pi^{*}}$. Then,
$T_{\pi}^{(t)}$ is a $t$-cone tree with $\pi$ as its degree sequence.
For example, let
$\pi_1=\left(18^{(2)},6^{(2)},5^{(4)},4,3^{(10)}\right)$ be a non-increasing degree sequence of a $2$-cone tree with $19$ vertices. Then,
$T_{\pi_1}^{(2)}$ is the corresponding  $2$-cone tree
as shown in Fig. 4.1.
\par\medskip
\[
\setlength{\unitlength}{0.5mm}
\begin{picture}(240,80)
\put(145,80){\circle*{3}}
\put(120,60){\circle*{3}}\put(170,60){\circle*{3}}
\put(65,60){\circle*{3}}\put(220,60){\circle*{3}}
\put(45,40){\circle*{3}}\put(65,40){\circle*{3}}\put(85,40){\circle*{3}}
\put(105,40){\circle*{3}}\put(135,40){\circle*{3}}
\put(155,40){\circle*{3}}\put(185,40){\circle*{3}}
\put(205,40){\circle*{3}}\put(235,40){\circle*{3}}
\put(35,20){\circle*{3}}\put(50,20){\circle*{3}}\put(65,20){\circle*{3}}

\put(146,81){$v_{0,1}$}
\put(123,57){$v_{1,2}$}\put(172,60){$v_{1,3}$}
\put(68,57){$v_{1,1}$}\put(222,60){$v_{1,4}$}
\put(46.5,35){$v_{2,1}$}\put(66.5,35){$v_{2,2}$}\put(86.5,35){$v_{2,3}$}
\put(106.5,35){$v_{2,4}$}\put(136.5,35){$v_{2,5}$}
\put(156.5,35){$v_{2,6}$}\put(186.5,35){$v_{2,7}$}
\put(206.5,35){$v_{2,8}$}\put(236.5,35){$v_{2,9}$}
\put(30,15){$v_{3,1}$}\put(45,15){$v_{3,2}$}\put(60,15){$v_{3,3}$}

\qbezier(145,80)(132.5,70)(120,60)\qbezier(145,80)(157.5,70)(170,60)
\qbezier(145,80)(105,70)(65,60)\qbezier(145,80)(172.5,70)(220,60)
\qbezier(120,60)(112.5,50)(105,40)\qbezier(120,60)(127.5,50)(135,40)
\qbezier(170,60)(162.5,50)(155,40)\qbezier(170,60)(177.5,50)(185,40)
\qbezier(65,60)(55,50)(45,40)\qbezier(65,60)(65,50)(65,40)\qbezier(65,60)(75,50)(85,40)
\qbezier(220,60)(212.5,50)(205,40)\qbezier(220,60)(227.5,50)(235,40)
\qbezier(45,40)(40,30)(35,20)\qbezier(45,40)(47.5,30)(50,20)\qbezier(65,40)(65,30)(65,20)

\put(0,65){\circle*{3}}\put(0,50){\circle*{3}}
\put(0,57.5){\circle{30}}\qbezier(0,50)(0,57.5)(0,65)\put(-11,56){$K_2$}
\put(0,65){\line(3,1){20}}\put(0,65){\line(1,0){20}}
\put(0,65){\line(4,3){20}}\put(22,68){$\vdots$}
\put(0,50){\line(3,-1){20}}\put(0,50){\line(1,0){20}}
\put(0,50){\line(4,-3){20}}\put(22,40){$\vdots$}

\put(145,80){\line(-3,2){10}}
\put(120,60){\line(-3,2){10}}\put(170,60){\line(-1,0){10}}
\put(65,60){\line(-3,2){10}}\put(220,60){\line(-1,0){10}}
\put(45,40){\line(-3,2){10}}\put(65,40){\line(-3,2){10}}\put(85,40){\line(-3,2){10}}
\put(105,40){\line(-3,2){10}}\put(135,40){\line(-3,2){10}}
\put(155,40){\line(-3,2){10}}\put(185,40){\line(-3,2){10}}
\put(205,40){\line(-3,2){10}}\put(235,40){\line(-3,2){10}}
\put(35,20){\line(-3,2){10}}\put(50,20){\line(-3,2){10}}\put(65,20){\line(-3,2){10}}

\put(145,80){\line(-3,1){10}}
\put(120,60){\line(-3,1){10}}\put(170,60){\line(-3,-1){10}}
\put(65,60){\line(-3,1){10}}\put(220,60){\line(-3,-1){10}}
\put(45,40){\line(-3,1){10}}\put(65,40){\line(-3,1){10}}\put(85,40){\line(-3,1){10}}
\put(105,40){\line(-3,1){10}}\put(135,40){\line(-3,1){10}}
\put(155,40){\line(-3,1){10}}\put(185,40){\line(-3,1){10}}
\put(205,40){\line(-3,1){10}}\put(235,40){\line(-3,1){10}}
\put(35,20){\line(-3,1){10}}\put(50,20){\line(-3,1){10}}\put(65,20){\line(-3,1){10}}

\put(90,0){Fig. 4.1 ~$T_{\pi_1}^{(2)}$}
\end{picture}
\]

\begin{thm}\label{42t}
If $c=0$, then  $G_{\pi}\cong T_{\pi}^{(t)}$.
\end{thm}
\begin{pf}
Let $G_{\pi}=K_t\vee T$. By Theorem \ref{35t} and $c=0$,
$T$ is a BFS-tree with tree degree sequence $\pi^{*}$.
As mentioned before, $T_{\pi^{*}}$ is the unique  BFS-tree with  tree  degree sequence   $\pi^{*}$.
So $T\cong T_{\pi^{*}}$,
and then $G_{\pi}\cong T_{\pi}^{(t)}$, as desired.
\end{pf}

In the following, we turn to consider the $\Theta_\alpha$-maximal unicyclic graph of $\mathscr{U}(\pi,t)$.
When
$\pi$ is the degree sequence  of a $t$-cone unicyclic graph,
then $\pi^{*}$ is a unicyclic graph degree sequence.
According to Proposition \ref{41p},
if $\pi$ satisfies $(4.1.1)$,
then $\pi=\left((n-1)^{(t)},(t+2)^{(n-t)}\right)$ and $\pi^*=\left(2^{(n-t)}\right)$,
in such case we let $U_{\pi}^{(t)}=K_t \vee C_{n-t}$;
if $\pi$ satisfies $(4.1.2)$,
then we first construct
a special unicyclic graph $U_{\pi^*}$ with $\pi^{*}$ as its degree sequence:

$(i)$ Put $s_0=1$, select a vertex $v_{0,1}$ as a root of $U_{\pi^*}$ and begin with $v_{0,1}$ of the zeroth layer;

$(ii)$ Put $s_1=d_{t+1}^*$, select $s_1$ vertices $\big\{v_{1,1},v_{1,2},\ldots,v_{1,s_1}\big\}$ of the first layer such that there are adjacent to $v_{0,1}$,
and $v_{1,1}$ is adjacent to $v_{1,2}$.
Then, $d_{U_{\pi^*}}(v_{0,1})=s_1=d_{t+1}^*$;

$(iii)$ Let $d_{U_{\pi^*}}(v_{1,i})=d_{t+i+1}^*$ for $i=1,2,\ldots,s_1$
\big(we shall stop if $t+i+1>n$\big).
Select $s_2=\sum_{i=1}^{s_1} d_{U_{\pi^*}}(v_{1,i})-s_1-2$ vertices $\big\{v_{2,1},v_{2,2},\ldots,v_{2,s_2}\big\}$
of the second layer such that $d_{U_{\pi^*}}(v_{1,1})-2$ vertices are adjacent to $v_{1,1}$,
$d_{U_{\pi^*}}(v_{1,2})-2$ vertices are adjacent to $v_{1,2}$ ,
and $d_{U_{\pi^*}}(v_{1,i})-1$ vertices are adjacent to $v_{1,i}$ for $i=3,4,\ldots,s_1$;

$(iv)$ The remaining vertices of $U_{\pi^*}$ appear in a BFS-connecting.

In this way, we obtain only one unicyclic graph $U_{\pi^*}$ with $\pi^*$ as its degree sequence.
It can be seen that $U_{\pi^*}$ is a BFS-unicyclic graph.
Besides,
$U_{\pi}^{(t)}:=K_t \vee U_{\pi^*}$ is a $t$-cone unicyclic graph with $\pi$ as its degree sequence.
As an illustrated example to $U_{\pi}^{(t)}$ for $\pi$ satisfying $(4.1.2)$, let
$\pi_3=\left(16^{(2)},6^{(2)},5^{(4)},4,3^{(8)}\right)$ be a non-increasing degree sequence of $2$-cone unicyclic graph with $17$ vertices, then  $U_{\pi_3}^{(2)}$ is the corresponding  $2$-cone unicyclic graph
as  shown in Fig. 4.2.
\par\medskip
\[
\setlength{\unitlength}{0.5mm}
\begin{picture}(240,80)
\put(145,80){\circle*{3}}
\put(120,60){\circle*{3}}\put(170,60){\circle*{3}}
\put(65,60){\circle*{3}}\put(220,60){\circle*{3}}
\put(45,40){\circle*{3}}\put(65,40){\circle*{3}}
\put(105,40){\circle*{3}}
\put(155,40){\circle*{3}}\put(185,40){\circle*{3}}
\put(205,40){\circle*{3}}\put(235,40){\circle*{3}}
\put(35,20){\circle*{3}}\put(50,20){\circle*{3}}\put(65,20){\circle*{3}}

\put(146,81){$v_{0,1}$}
\put(123,57){$v_{1,2}$}\put(172,60){$v_{1,3}$}
\put(68,55){$v_{1,1}$}\put(222,60){$v_{1,4}$}
\put(46.5,35){$v_{2,1}$}\put(66.5,35){$v_{2,2}$}
\put(106.5,35){$v_{2,3}$}
\put(156.5,35){$v_{2,4}$}\put(186.5,35){$v_{2,5}$}
\put(206.5,35){$v_{2,6}$}\put(236.5,35){$v_{2,7}$}
\put(30,15){$v_{3,1}$}\put(45,15){$v_{3,2}$}\put(60,15){$v_{3,3}$}

\qbezier(145,80)(132.5,70)(120,60)\qbezier(145,80)(157.5,70)(170,60)
\qbezier(145,80)(105,70)(65,60)\qbezier(145,80)(172.5,70)(220,60)
\qbezier(120,60)(112.5,50)(105,40)
\qbezier(170,60)(162.5,50)(155,40)\qbezier(170,60)(177.5,50)(185,40)
\qbezier(65,60)(55,50)(45,40)\qbezier(65,60)(65,50)(65,40)
\qbezier(220,60)(212.5,50)(205,40)\qbezier(220,60)(227.5,50)(235,40)
\qbezier(45,40)(40,30)(35,20)\qbezier(45,40)(47.5,30)(50,20)\qbezier(65,40)(65,30)(65,20)

\put(0,65){\circle*{3}}\put(0,50){\circle*{3}}
\put(0,57.5){\circle{30}}\qbezier(0,50)(0,57.5)(0,65)\put(-11,56){$K_2$}
\put(0,65){\line(3,1){20}}\put(0,65){\line(1,0){20}}
\put(0,65){\line(4,3){20}}\put(22,68){$\vdots$}
\put(0,50){\line(3,-1){20}}\put(0,50){\line(1,0){20}}
\put(0,50){\line(4,-3){20}}\put(22,40){$\vdots$}

\put(145,80){\line(-3,2){10}}
\put(120,60){\line(-3,2){10}}\put(170,60){\line(-1,0){10}}
\put(65,60){\line(-3,2){10}}\put(220,60){\line(-1,0){10}}
\put(45,40){\line(-3,2){10}}\put(65,40){\line(-3,2){10}}
\put(105,40){\line(-3,2){10}}
\put(155,40){\line(-3,2){10}}\put(185,40){\line(-3,2){10}}
\put(205,40){\line(-3,2){10}}\put(235,40){\line(-3,2){10}}
\put(35,20){\line(-3,2){10}}\put(50,20){\line(-3,2){10}}\put(65,20){\line(-3,2){10}}

\put(145,80){\line(-3,1){10}}
\put(120,60){\line(-3,1){10}}\put(170,60){\line(-3,-1){10}}
\put(65,60){\line(-3,1){10}}\put(220,60){\line(-3,-1){10}}
\put(45,40){\line(-3,1){10}}\put(65,40){\line(-3,1){10}}
\put(105,40){\line(-3,1){10}}
\put(155,40){\line(-3,1){10}}\put(185,40){\line(-3,1){10}}
\put(205,40){\line(-3,1){10}}\put(235,40){\line(-3,1){10}}
\put(35,20){\line(-3,1){10}}\put(50,20){\line(-3,1){10}}\put(65,20){\line(-3,1){10}}

\qbezier(65,60)(80,60)(120,60)

\put(90,0){Fig. 4.2 ~$U_{\pi_3}^{(2)}$}
\end{picture}
\]

\begin{prop}\label{43p}
Suppose that $\pi^*$ satisfies $(4.1.2)$,
and let $U^*$ be a BFS-unicyclic graph with $\pi^*$ as its degree sequence such that $\mathcal{B}(U^*)=C_3$ with $V(C_3)=\big\{v_{t+1},v_{t+2},v_{t+3}\big\}$.
Then $U^*\cong U_{\pi^*}$.
\end{prop}

\begin{pf}
Note that $d_{t+1}^*\geq 3$ and $d_{t+2}^*\geq 2$ by $(4.1.2)$.
If $d_{t+1}^*=3$ and $d_{t+2}^*=2$,
since $\mathcal{B}(U^*)=C_3$,
then $U^*$ shall be obtained from $P_{n-t-3}$ and $C_3$ by adding one edge between one end vertex of $P_{n-t-3}$ and one vertex of $C_3$, and it can be checked that $U^*\cong U_{\pi^*}$.
Otherwise,
$d_{t+1}^*\geq 4$ or $d_{t+2}^*\geq 3$.
Let $\pi_2=\left(d_{t+1}^*+d_{t+2}^*+d_{t+3}^*-6,d_{t+4}^*,
d_{t+5}^*,\ldots,d_{n}^*\right)
:=\left(d''_1,d''_2,d''_3,\ldots,d''_{n-t-2}\right)$.
Since
$\sum_{i=1}^{n-t-2}d''_i=\sum_{j=t+1}^{n}d_{j}^*-6=2\big((n-t-2)-1\big)$
and $d''_1=d_{t+1}^*+d_{t+2}^*+d_{t+3}^*-6\geq d_{t+3}^*\geq d''_2=d_{t+4}^*\geq d''_3=d_{t+5}^*\geq\cdots\geq d''_{n-t-2}=d_{n}^*$,
then $\pi_2$ is a non-increasing tree degree sequence and so there is a unique BFS-tree $T_{\pi_2}$
with $\pi_2$ as its degree sequence by the former arguments.
Clearly,
$U^*$ and $U_{\pi^*}$ are both BFS-unicyclic graphs with $\pi^*$ as their degree sequence such that $\mathcal{B}(U^*)=\mathcal{B}(U_{\pi^*})=C_3$ with $V(C_3)=\big\{v_{t+1},v_{t+2},v_{t+3}\big\}$.
So the graph obtained from $U^*$ \big(resp., $U_{\pi^*}$\big) by contracting the $C_3$ is isomorphic to $T_{\pi_2}$.
Now, it is easy to see that the positions of $\big\{v_{t+4},v_{t+5},\ldots,v_n\big\}$ of $U^*$ \big(resp., $U_{\pi^*}$\big) is fixed by the uniqueness of $T_{\pi_2}$,
and so $U^*\cong U_{\pi^*}$.
\end{pf}

\begin{thm}\label{44t}
If $c=1$ and  $0\leq t\leq n-3$,
then $G_{\pi}\cong U_{\pi}^{(t)}$.
\end{thm}

\begin{pf}
If $\pi$ satisfies $(4.1.1)$,
then $\mathscr{U}(\pi,t)=\{K_t \vee C_{n-t}\}$,
and hence $G_{\pi}\cong U_{\pi}^{(t)}=K_t \vee C_{n-t}$, as desired.
Otherwise, $\pi$ satisfies $(4.1.2)$.
Let $G_{\pi}=K_t\vee H$.
By Theorem \ref{35t}, we may suppose that $v_{t+1}\prec v_{t+2}\prec\cdots \prec v_n$ is a   good   BFS-ordering of $V(H)$.

Since $c=1$, we  let   $C_g$ be the unique cycle of $H$,
where $g\geq 3$.
Notice that $\{v_{t+1},v_{t+2},v_{t+3}\}\subseteq V(C_g)$ by Theorem \ref{37t}.
Next, we need to show that $g=3$.
By contradiction, suppose that $g\geq 4$,
then there exists a vertex $x\in V(C_g)\setminus N_H(v_{t+1})$.
Since $C_g$ is the unique cycle of $H$ and
since $d^{*}_{t+1}\geq 3$ by $(4.1.2)$,
  there exists a vertex $y\in N_H(v_{t+1})\setminus V(C_g)$.
On one hand,
since $y\in N_H(v_{t+1})$ and $x\notin N_H(v_{t+1})$,
then $h_H(y)<h_H(x)$,
which leads to $f_{G_{\pi}}(y)\geq f_{G_{\pi}}(x)$ by Definition \ref{26d} $(ii)$.
On the other hand,
since $x\in V(C_g)$ and $y\notin V(C_g)$,
then by Theorem \ref{37t},
we get $f_{G_{\pi}}(x)>f_{G_{\pi}}(y)$,  a contradiction.
Therefore,
we conclude that $g=3$ and then
$V(C_3)=\big\{v_{t+1},v_{t+2},v_{t+3}\big\}$.
Now by Theorem \ref{35t},
$H$ is a BFS-unicyclic graph with $\mathcal{B}(H)\cong C_3$ and $V(C_3)=\big\{v_{t+1},v_{t+2},v_{t+3}\big\}$,
so $H\cong U_{\pi^*}$ by Proposition  \ref{43p},
and it follows that ${G_{\pi}}\cong U_{\pi}^{(t)}$, as desired.
\end{pf}

Finally, we shall discuss the $\Theta_\alpha$-maximal bicyclic graph of $\mathscr{B}(\pi,t)$.
We first introduce three notations and we will use them hereafter without specified indicated:
let $B(n_1,n_2)$ denote the bicyclic  graph obtained from two cycles $C_{n_1}$ and $C_{n_2}$
by adding one edge  connecting a vertex of $C_{n_1}$ and  a vertex of $C_{n_2}$;
let $C(n_1,n_2)$ denote the bicyclic  graph obtained from two cycles $C_{n_1}$ and $C_{n_2}$
by identifying a vertex of $C_{n_1}$ with a vertex of $C_{n_2}$;
and let $\theta(p,r,q)$ denote the bicyclic  graph obtained from three vertex-disjoint paths, say
 $P_{p+1}$, $P_{r+1}$, and  $P_{q+1}$, respectively,
by identifying the three initial (resp. terminal) vertices of them, where $q\geq r\geq 1$ and $p\geq q\geq 2$.

When $\pi$ is the degree sequence
of a $t$-cone bicyclic graph,
then $\pi^{*}$ is a bicyclic graph degree sequence.
By Proposition \ref{41p},
$\pi$ should be one of the following four cases,
and we can construct a special $t$-cone bicyclic graph $B_{\pi}^{(t)}$ with $\pi$ as its degree sequence:

\noindent\textbf{Case 1.~} $\pi$ satisfies $(4.2.1)$.
Let $B_{\pi}^{(t)}=K_t\vee C(3,n-t-2)$;

\noindent\textbf{Case 2.~} $\pi$ satisfies $(4.2.2)$.
Let $B_{\pi}^{(t)}=K_t\vee \theta(n-t-2,1,2)$;

\noindent\textbf{Case 3.~} $\pi$ satisfies $(4.2.3)$.
Let $B_{\pi}^{(t)}=K_t\vee B^\star$,
where $B^\star$ is the bicyclic graph with $n-t$ vertices obtained from $C(3,3)$ by attaching $d_{t+1}^*-4$ paths of almost equal lengths \big(i.e. their lengths are different at most by one\big) to the maximum degree vertex of $C(3,3)$;

\noindent\textbf{Case 4.~} $\pi$ satisfies $(4.2.4)$.
Then we shall construct a $t$-cone bicyclic graph $B_{\pi}^{(t)}=K_t\vee B_{\pi^*}$ such that $\pi^{*}$ is the degree sequence of $B_{\pi^*}$:

$(i)$ Put $s_0=1$, select a vertex $v_{0,1}$ as a root of $B_{\pi^*}$ and begin with $v_{0,1}$ of the zeroth layer;

$(ii)$ Put $s_1=d_{t+1}^*$, select $s_1$ vertices $\big\{v_{1,1},v_{1,2},\ldots,v_{1,s_1}\big\}$ of the first layer such that there are adjacent to $v_{0,1}$,
 and  $v_{1,1}$ is adjacent to $v_{1,2}$ and  $v_{1,3}$.
Then, $d_{B_{\pi^*}}(v_{0,1})=s_1=d_{t+1}^*$;

$(iii)$ Let $d_{B_{\pi^*}}(v_{1,i})=d_{t+i+1}^*$ for $i=1,2,\ldots,s_1$
\big(we shall stop if $t+i+1>n$\big).
Select $s_2=\sum_{i=1}^{s_1} d_{B_{\pi^*}}(v_{1,i})-s_1-4$ vertices $\big\{v_{2,1},v_{2,2},\ldots,v_{2,s_2}\big\}$
of the second layer such that $d_{B_{\pi^*}}(v_{1,1})-3$ vertices are adjacent to $v_{1,1}$,
$d_{B_{\pi^*}}(v_{1,2})-2$ vertices are adjacent to $v_{1,2}$,
$d_{B_{\pi^*}}(v_{1,3})-2$ vertices are adjacent to $v_{1,3}$,
and $d_{B_{\pi^*}}(v_{1,i})-1$ vertices are adjacent to $v_{1,i}$ for $i=4,5,\ldots,s_1$;

$(iv)$ The remaining vertices of $B_{\pi^*}$ appear in a BFS-connecting.

In this way,
we obtain the unique  bicyclic graph $B_{\pi^*}$ with $\pi^*$ as its degree sequence.
It can be seen that $B_{\pi^*}$ is a BFS-bicyclic graph.
Moreover,
$B_{\pi}^{(t)}=K_t \vee B_{\pi^*}$ is a $t$-cone bicyclic graph with $\pi$ as its degree sequence.
For instance,
$\pi_4=\left(14^{(2)},6^{(2)},5^{(4)},4,3^{(6)}\right)$ is a non-increasing degree sequence of $2$-cone bicyclic graph with $15$ vertices,
and $B_{\pi_4}^{(2)}$ is the corresponding $2$-cone bicyclic graph
as  shown in Fig. 4.3.
\par\medskip
\[
\setlength{\unitlength}{0.5mm}
\begin{picture}(240,80)
\put(145,80){\circle*{3}}
\put(120,60){\circle*{3}}\put(170,60){\circle*{3}}
\put(65,60){\circle*{3}}\put(220,60){\circle*{3}}
\put(45,40){\circle*{3}}
\put(105,40){\circle*{3}}
\put(155,40){\circle*{3}}
\put(205,40){\circle*{3}}\put(235,40){\circle*{3}}
\put(35,20){\circle*{3}}\put(50,20){\circle*{3}}\put(105,20){\circle*{3}}

\put(146,81){$v_{0,1}$}
\put(123,57){$v_{1,2}$}\put(172,60){$v_{1,3}$}
\put(68,53){$v_{1,1}$}\put(222,60){$v_{1,4}$}
\put(46.5,35){$v_{2,1}$}
\put(106.5,35){$v_{2,2}$}
\put(156.5,35){$v_{2,3}$}
\put(206.5,35){$v_{2,4}$}\put(236.5,35){$v_{2,5}$}
\put(30,15){$v_{3,1}$}\put(45,15){$v_{3,2}$}\put(106.5,15){$v_{3,3}$}

\qbezier(145,80)(132.5,70)(120,60)\qbezier(145,80)(157.5,70)(170,60)
\qbezier(145,80)(105,70)(65,60)\qbezier(145,80)(172.5,70)(220,60)
\qbezier(120,60)(112.5,50)(105,40)
\qbezier(170,60)(162.5,50)(155,40)
\qbezier(65,60)(55,50)(45,40)
\qbezier(220,60)(212.5,50)(205,40)\qbezier(220,60)(227.5,50)(235,40)
\qbezier(45,40)(40,30)(35,20)\qbezier(45,40)(47.5,30)(50,20)\qbezier(105,40)(105,30)(105,20)

\put(0,65){\circle*{3}}\put(0,50){\circle*{3}}
\put(0,57.5){\circle{30}}\qbezier(0,50)(0,57.5)(0,65)\put(-11,56){$K_2$}
\put(0,65){\line(3,1){20}}\put(0,65){\line(1,0){20}}
\put(0,65){\line(4,3){20}}\put(22,68){$\vdots$}
\put(0,50){\line(3,-1){20}}\put(0,50){\line(1,0){20}}
\put(0,50){\line(4,-3){20}}\put(22,40){$\vdots$}

\put(145,80){\line(-3,2){10}}
\put(120,60){\line(-1,0){12}}\put(170,60){\line(-1,0){10}}
\put(65,60){\line(-3,2){10}}\put(220,60){\line(-1,0){10}}
\put(45,40){\line(-3,2){10}}
\put(105,40){\line(-3,2){10}}
\put(155,40){\line(-3,2){10}}
\put(205,40){\line(-3,2){10}}\put(235,40){\line(-3,2){10}}
\put(35,20){\line(-3,2){10}}\put(50,20){\line(-3,2){10}}\put(105,20){\line(-3,2){10}}

\put(145,80){\line(-3,1){10}}
\put(120,60){\line(-3,1){10}}\put(170,60){\line(-3,-1){10}}
\put(65,60){\line(-3,1){10}}\put(220,60){\line(-3,-1){10}}
\put(45,40){\line(-3,1){10}}
\put(105,40){\line(-3,1){10}}
\put(155,40){\line(-3,1){10}}
\put(205,40){\line(-3,1){10}}\put(235,40){\line(-3,1){10}}
\put(35,20){\line(-3,1){10}}\put(50,20){\line(-3,1){10}}\put(105,20){\line(-3,1){10}}

\qbezier(65,60)(95,50)(120,60)
\qbezier(65,60)(120,71)(170,60)

\put(90,-5){Fig. 4.3 ~$B_{\pi_4}^{(2)}$}
\end{picture}
\]
\par\medskip

Analogously as the proof of Proposition \ref{43p},
we have the following statement.

\begin{prop}\label{45p}
Suppose that $\pi^*$ satisfies $(4.2.4)$,
and let $B^*$ be a BFS-bicyclic graph with $\pi^*$ as its degree sequence such that $\mathcal{B}(B^*)=\theta(2,1,2)$ with $V\big(B^*\big)=\big\{v_{t+1},v_{t+2},v_{t+3},v_{t+4}\big\}$, $v_{t+1}$ and $v_{t+2}$ are the two vertices of degree three in $\theta(2,1,2)$.
Then, $B^*\cong B_{\pi^*}$.
\end{prop}

\begin{thm}\label{46t}
If $c=2$ and $0\leq t\leq n-4$, then $G_{\pi}\cong B_{\pi}^{(t)}$.
\end{thm}

\begin{pf}
Let $G_{\pi}=K_t\vee H$.
By Proposition \ref{41p},
we need to consider the following four cases.

\noindent\textbf{Case 1.~} $\pi$ satisfies $(4.2.1)$.
Then $H\cong C(p,n-t+1-p)$, where $n-t+1-p\geq p\geq 3$.
By Theorem \ref{35t},
we may assume that
$N_H(v_{t+1})=\big\{w_{1},w_{2},w_{3},w_{4}\big\}$,
$\big\{w_{1},w_{2}\big\}\subseteq V(C_p)$ and
$\big\{w_{3},w_{4}\big\}\subseteq V(C_{n-t+1-p})$.
Clearly,
$d_H(v_{t+1})=4>d_H(w_{1})=d_H(w_{2})=d_H(w_{3})=d_H(w_{4})=2$.
Suppose that $H\ncong C(3,n-t-2)$.
Then there exist vertices $x$ and $y$ such that
$x\in N_H(w_{1})\setminus \{v_{t+1},w_{2}\}$
and $y\in N_H(w_{3})\setminus \{v_{t+1},w_{4}\}$.
Let $H_1=H+v_{t+1}x+w_{1}w_{3}-w_{1}x-v_{t+1}w_{3}$ and $G_1=K_t\vee H_1$.
Since $f_{G_{\pi}}(v_{t+1})>f_{G_{\pi}}(w_{1})$ by $(\ref{31e})$ and Lemma \ref{31l},
then by the fact that $H_1$ is connected and by Corollary \ref{24c},
$f_{G_{\pi}}(x)<f_{G_{\pi}}(w_{3})$.
Let $H_2=H+xy+w_{1}w_{3}-w_{3}y-w_{1}x$ and $G_2=K_t\vee H_2$.
Clearly,
$G_2=G+xy+w_{1}w_{3}-w_{3}y-w_{1}x\in \mathscr{B}(\pi,t)$.
Since $f_{G_{\pi}}(w_{1})\geq f_{G_{\pi}}(y)$ by Theorem \ref{35t} as $h_H(w_1)<h_H(y)$, and recall that $f_{G_{\pi}}(w_{3})>f_{G_{\pi}}(x)$, it follows from Theorem \ref{22t} that $\Theta(G_2,\alpha)>\Theta({G_{\pi}},\alpha)$,
which is contrary with the definition of $G_{\pi}$.
Thus,
$H\cong C(3,n-t-2)$ and so $G_{\pi}\cong B_{\pi}^{(t)}=K_t\vee C(3,n-t-2)$.

\noindent\textbf{Case 2.~} $\pi$ satisfies $(4.2.2)$.
Suppose that $H\ncong \theta(n-t-2,1,2)$.
Then by Theorem \ref{35t},
either $H\cong \theta(n-t-p,1,p)$ or $H\cong B(p,n-p-t)$,
where $n-p-t\geq p\geq 3$,
$N_H(v_{t+1})=\{v_{t+2},v_{t+3},v_{t+4}\}$,
and $d_H(v_{t+1})=d_H(v_{t+2})=3>d_H(v_{t+3})=d_H(v_{t+4})=2$.
Bearing in mind that $H\ncong \theta(n-t-2,1,2)$,
then there exist vertices $x$ and $y$ such that $x\in N_H(v_{t+3})\setminus \big\{v_{t+1},v_{t+2}\big\}$ and $y\in N_H(v_{t+2})\setminus \big(\big\{v_{t+1},v_{t+3},v_{t+4},x\big\}\cup N_H(x)\big)$.
Let $H_3=H+v_{t+2}v_{t+3}+xy-v_{t+2}y-v_{t+3}x$ and $G_3=K_t\vee H_3$,
whenever $H\cong \theta(1,p,n-t-p)$ or $H\cong B(p,n-p-t)$.
Then,  $H_3$ is connected
and so $G_3\in \mathscr{B}(\pi,t)$.
By $(\ref{31e})$, Lemma \ref{31l} and Theorem \ref{35t},
we have $f_{G_{\pi}}(v_{t+2})>f_{G_{\pi}}(x)$ since $d_H(v_{t+2})=3>d_H(x)=2$,
and $f_{G_{\pi}}(v_{t+3})\geq f_{G_{\pi}}(y)$ since $y\notin N_H(v_{t+1})=\big\{v_{t+1},v_{t+3},v_{t+4}\big\}$.
Thus,
Theorem \ref{22t} implies that $\Theta(G_3,\alpha)>\Theta({G_{\pi}},\alpha)$,
contrary with  the choice of ${G_{\pi}}$.
This contradiction confirms that  ${G_{\pi}}\cong B_{\pi}^{(t)}=K_t\vee \theta(n-t-2,1,2)$.

\noindent\textbf{Case 3.~} $\pi$ satisfies $(4.2.3)$.
Since $d_n^*=1$,
it follows from  Theorem \ref{310t} that  $\mathcal{B}(H)\cong C(3,3)$.
Combining this with Theorem \ref{35t},
it is easy to see that $G_{\pi}\cong B_{\pi}^{(t)}=K_t\vee B^\star$,
where $B^\star$ is the bicyclic graph with $n-t$ vertices obtained from $C(3,3)$ by attaching $d_{t+1}^*-4$ paths of almost equal lengths to the maximum degree vertex of $C(3,3)$.

\noindent\textbf{Case 4.~} $\pi$ satisfies $(4.2.4)$.
According to Theorems \ref{35t} and \ref{311t},
we just need to consider the following three subcases:

\noindent\textbf{Subcase 4.1~} $\mathcal{B}(H)\cong B(3,3)$.
By Theorems \ref{37t} and \ref{38t},
we have
$N_H\big(v_{t+1}\big)\cap V\big(B(3,3)\big)=\big\{v_{t+2},v_{t+3},v_{t+4}\big\}$
and $d_{B(3,3)}\big(v_{t+1}\big)=d_{B(3,3)}\big(v_{t+2}\big)=3>d_{B(3,3)}\big(v_{t+3}\big)=d_{B(3,3)}\big(v_{t+4}\big)=2$.
Suppose that $x\in N_{B(3,3)}\big(v_{t+2}\big)$ and $x\neq v_{t+1}$.
Let $H_4=H+v_{t+2}v_{t+3}+v_{t+4}x-v_{t+3}v_{t+4}-v_{t+2}x$ and $G_4=K_t\vee H_4$.
Since $H_4$ is connected,
$G_4=K_t\vee H_4\in \mathscr{B}(\pi,t)$.
Note that $f_{G_{\pi}}\big(v_{t+2}\big)>f_{G_{\pi}}\big(v_{t+4}\big)$ by Theorem \ref{38t},
and $f_{G_{\pi}}\big(v_{t+3}\big)\geq f_{G_{\pi}}(x)$ by Theorem \ref{35t}
\big(since $v_{t+3}\in N_H(v_{t+1})$ and $x\notin N_H\big(v_{t+1}\big)$\big).
Hence we get $\Theta(G_4,\alpha)>\Theta(G_{\pi},\alpha)$ by Theorem \ref{22t},
contrary with  the definition  of ${G_{\pi}}$.

\noindent\textbf{Subcase 4.2~} $\mathcal{B}(H)\cong C(3,3)$.
It follows from Theorems \ref{37t} and \ref{38t} that
$N_H(v_{t+1})\cap V\big(C(3,3)\big)=\big\{v_{t+2},v_{t+3},v_{t+4},v_{t+5}\big\}$
and $d_{C(3,3)}\big(v_{t+1}\big)=4>d_{C(3,3)}\big(v_{t+2}\big)=d_{C(3,3)}\big(v_{t+3}\big)=d_{C(3,3)}\big(v_{t+4}\big)=d_{C(3,3)}\big(v_{t+5}\big)=2$.
Since $d_H\big(v_{t+2}\big)\geq 3$ by Theorem \ref{35t} and $(4.2.4)$,
there exists a vertex $u$ such that $u\in N_H\big(v_{t+2}\big)\setminus V\big(C(3,3)\big)$.
We may suppose that $\big\{v_{t+5},v_{t+4}\big\}\not\in N(v_{t+2})$.
Let $H_5=H+v_{t+2}v_{t+4}+v_{t+5}u-v_{t+4}v_{t+5}-uv_{t+2}$ and $G_5=K_t\vee H_5$.
It can be easily checked  that $H_5$ is connected,
and thus $G_5\in \mathscr{B}(\pi,t)$.
Since $f_{G_{\pi}}(v_{t+4})>f_{G_{\pi}}(u)$ by Theorem \ref{37t},
and $f_{G_{\pi}}(v_{t+2})\geq f_{G_{\pi}}(v_{t+5})$ by Theorem \ref{35t},
we have $\Theta(G_5,\alpha)>\Theta({G_{\pi}},\alpha)$ by Theorem \ref{22t},
also contrary with  the choice of ${G_{\pi}}$.

\noindent\textbf{Subcase 4.3~} $\mathcal{B}(H)\cong \theta(2,1,2)$.
By Theorems \ref{35t}, \ref{37t} and \ref{38t},
$H$ is a BFS-bicyclic graph
such that $\mathcal{B}(H)=\theta(2,1,2)$, $V\big(H\big)=\big\{v_{t+1},v_{t+2},v_{t+3},v_{t+4}\big\}$, $v_{t+1}$ and $v_{t+2}$ are the two vertices of degree three in $\theta(2,1,2)$.
Hence $H\cong B_{\pi^*}$ by Proposition \ref{45p},
and so $G_{\pi}\cong B_{\pi}^{(t)}=K_t\vee B_{\pi^*}$.
\end{pf}


\section{\large{Majorization theorems for the general spectral radius of $t$-cone trees, $t$-cone unicyclic graphs and $t$-cone bicyclic graphs}}

If ${\bf y}=(y_1,y_2,\ldots,y_n)$ is a  non-increasing integer sequence and $y_i\geq y_j+2$, then the following operation is called a {\em unit transformation} from $i$ to $j$ on ${\bf y}$: subtract 1 from $y_i$ and add 1 to $y_j$. To begin with,
we need to introduce the following important result.

\begin{lem}\label{51l}{\em\cite{Marshall}}
Let
${\bf y}=\big(y_1,y_2,\ldots,y_n\big)$ and ${\bf z}=\big(z_1,z_2,\ldots,z_n\big)$   be two positive non-increasing integer sequences. If ${\bf y}\lhd {\bf z}$, then   ${\bf y}$ can be obtained from ${\bf z}$ by a finite sequence of unit transformations.
\end{lem}

Recall that   $\pi=\big(d_1,d_2,\ldots,d_n\big)$ and $\pi'=\big(d'_1, d'_2,\ldots,d'_n\big)$ always denote two non-increasing degree sequences of $t$-cone $c$-cyclic graphs  such that  $\pi\lhd \pi'$, and
$G_{\pi}$ and $G_{\pi'}$ are  the $\Theta_\alpha$-maximal graphs of $\Gamma(\pi,t;c)$ and $\Gamma(\pi',t;c)$, respectively.
If $\pi$ and $\pi'$ differ only in two positions where the difference is $1$
\big(namely, $\pi$ can be obtained from $\pi'$ by a unit transformation\big),
then $\pi\lhd \pi'$ is called a {\em star majorization} and denoted by  $\pi\lhd^* \pi'$.
This implies that,
if $\pi\lhd^* \pi'$,
then  $d_p'=d_p+1$,
$d_q'=d_q-1$,
and $d_i=d_i'$ for $i\in\big\{1,2,\ldots,n\big\}\setminus \big\{p,q\big\}$,
where $1\leq p<q\leq n$.

\begin{rem}\label{52r}
{\em By Lemma \ref{51l}, if   $\pi\lhd \pi'$, then there exist a series   of  non-increasing degree  sequences of $t$-cone $c$-cyclic graphs,
say $\pi_1,\pi_2,\ldots,\pi_{k-1}$,
such that $(\pi=)$ $\pi_0\lhd\pi_1\lhd\cdots\lhd\pi_{k-1}\lhd\pi_k$ $(=\pi')$,
and $\pi_i\lhd^*\pi_{i+1}$ holds for any $i\in \big\{0,1,\ldots,k-1\big\}$.
Therefore,
without loss of generality,
we may just simplify  assume that $\pi\lhd^* \pi'$ if $\pi\lhd \pi'$ when  considering the majorization theorems in the sequel. }
\end{rem}

\begin{rem}\label{53r}
{\em
If $\pi\lhd^* \pi'$,
since $d_j=d'_j=n-1$ for $1\leq j\leq t$,
we may suppose that
$d_p'=d_p+1$,
$d_q'=d_q-1$,
and $d_i=d_i'$ for $i\in \big\{1,2,\ldots, n\big\}\setminus \big\{p,q\big\}$,
where $t+1\leq p<q\leq n$.
Let $G_{\pi}=K_t\vee H$ and recall that $P_H\big(v_p,v_q\big)$ is a shortest path from $v_p$ to $v_q$.
By Theorem \ref{35t},
there is a good BFS-ordering $v_{t+1}\prec v_{t+2}\prec\cdots\prec v_n$ of $V(H)$ such that $d_{G_{\pi}}(v_i)=d_i$ for $i\in \big\{t+1,t+2,\ldots, n\big\}$.
If there is a vertex $w\in V(H)$ such that $w\in N_H(v_q)\setminus \big(N_H(v_p)\cup \{v_p\}\big)$ and $w\notin V\big(P_H(v_p,v_q)\big)$,
then we call $w$ a  {\em surprising vertex} of $G_{\pi}$. }
\end{rem}

\begin{lem}\label{54l}
If $\pi\lhd^* \pi'$ and  $G_{\pi}$ contains a surprising vertex,
then $\Theta(G_{\pi},\alpha)<\Theta(G_{\pi'},\alpha)$.
\end{lem}

\begin{pf}
We follow the notations in Remark \ref{53r}.
If $G_{\pi}=K_t\vee H$ contains a surprising vertex, say $w$, then
let $H_1=H+v_pw-v_qw$ and $G_1=K_t\vee H_1$.
Since $w\notin V\big(P_H(v_p,v_q)\big)$,
$H_1$ is connected,
and so $G_1\in \Gamma(\pi',t;c)$.
It follows from Theorem \ref{35t} that $f_{G_{\pi}}(v_p)\geq f_{G_{\pi}}(v_q)$, as $p<q$.
By Corollary \ref{23c} and the definition  of $G_{\pi'}$,
$\Theta({G_{\pi}},\alpha)<\Theta(G_1,\alpha)\leq\Theta(G_{\pi'},\alpha)$.
\end{pf}

When $c\in \big\{0,1,2\big\}$,  it is easily checked that there is at most  one non-increasing degree sequence of $t$-cone $c$-cyclic graph when $n-t\leq c+2$, and so
we may suppose that $n-t\geq c+3$ in the sequel.
Now we are ready to  state our main results of this section.

\begin{thm}\label{55t}
If $c=0$ and $0\leq t\leq n-3$, then   $\Theta(G_{\pi},\alpha)<\Theta(G_{\pi'},\alpha)$.
\end{thm}

\begin{pf}
Since $\pi\lhd \pi'$,
we may assume that $\pi\lhd^* \pi'$ by   Remark \ref{52r}.
Following the notations in Remark \ref{53r},
and according to $(4.0.1)$,
we can conclude that
$d_p'=d_p+1$,
$d_q=d_q'+1\geq t+2$, and
$d_j=d_j'\geq t+1$ for $j\in \big\{t+1,t+2,\ldots,n\big\}\setminus \big\{p,q\big\}$,
where $t+1\leq p<q\leq n$.
By Theorem \ref{42t},
we have $G_{\pi}\cong K_t \vee T_{\pi^{*}}$,
where $T_{\pi^{*}}$ is the unique   BFS-tree with $\pi^{*}$ as its degree sequence.
Note that $N_{T_{\pi^{*}}}\big(v_q\big)$ contains exactly one vertex in the path   $P_{T_{\pi^{*}}}\big(v_p,v_q\big)$.
Combining this with $d_{T_{\pi^{*}}}(v_q)=d_{G}(v_q)-t=d_q-t\geq 2$,
it is easily checked that  $G_{\pi}$ contains a surprising vertex, and hence  $\Theta(G_{\pi},\alpha)<\Theta(G_{\pi'},\alpha)$ by Lemma \ref{54l}.
\end{pf}

\begin{thm}\label{56t}
If $c=1$ and $0\leq t\leq n-4$, then  $\Theta(G_{\pi},\alpha)<\Theta(G_{\pi'},\alpha)$.
\end{thm}
\begin{pf}
Since $\pi\lhd \pi'$,
we may suppose that $\pi\lhd^* \pi'$ by Remark \ref{52r}.

\noindent\textbf{Case 1.~} $\pi$ satisfies $(4.1.1)$.
In this case,
$$\pi=\left((n-1)^{(t)}, (t+2)^{(n-t)}\right)~~\text{and} ~~ \pi'=\left((n-1)^{(t)},
t+3,(t+2)^{(n-t-2)},t+1\right).$$
So $G_{\pi}=K_t\vee C_{n-t}$.
Since $n-t\geq 4$,
denote by $C_{n-t}=w_1w_2w_3w_4\cdots w_{n-t}w_1$.
If $f_{G_{\pi}}(w_2)\geq f_{G_{\pi}}(w_3)$,
then let $H_1=C_{n-t}-w_4w_3+w_4w_2$ and $G_1=K_t\vee H_1$.
It can be seen that $G_1\in \mathscr{U}(\pi',t)$.
By Corollary \ref{23c},
we get $\Theta(G_{\pi},\alpha)<\Theta(G_1,\alpha)\leq\Theta(G_{\pi'},\alpha)$.
If $f_{G_{\pi}}(w_2)\leq f_{G_{\pi}}(w_3)$,
then let $H_2=C_{n-t}-w_1w_2+w_1w_3$ and $G_2=K_t\vee H_2$.
Note that
$G_2\in \mathscr{U}(\pi',t)$.
Once again,
Corollary \ref{23c} implies that  $\Theta(G_{\pi},\alpha)<\Theta(G_2,\alpha)\leq\Theta(G_{\pi'},\alpha)$.

\noindent\textbf{Case 2~} $\pi$ satisfies $(4.1.2)$.
In this case, by Theorem \ref{44t} we have  $G_{\pi}\cong K_t\vee U_{\pi^*}$,
where $U_{\pi^*}$ is a unique   BFS-unicyclic graph with  $\pi^*$ as its degree sequence  and $\mathcal{B}(U_{\pi^*})=C_3$.
Following the notations in Remark \ref{53r},
we may suppose that $V(C_3)=\big\{v_{t+1},v_{t+2},v_{t+3}\big\}$,
$d_p'=d_p+1$,
$d_q=d_q'+1$, and
$d_i=d_i'$ for $i\in \big\{t+1,t+2,\ldots,n\big\}\setminus \big\{p,q\big\}$,
where $t+1\leq p<q\leq n$.

\noindent\textbf{Subcase 2.1~} $q=t+2$ or $t+3$.
Note that $d_{t+1}'\geq d_{t+2}'\geq d_{t+3}'\geq t+2$, as $G_{\pi'}=K_t\vee H'$ and $H'$ contains a cycle.
Thus, $d_q=d_q'+1\geq t+3$  and so  $d_{U_{\pi^*}}\big(v_q\big)=d_q-t\geq 3$,
we have $v_p\in V(C_3)$ and $N_{U_{\pi^*}}\big(v_q\big)\setminus V(C_3)\neq \emptyset$.
Thus,
there is a vertex $w\in N_{U_{\pi^*}}\big(v_q\big)\setminus V(C_3)$ such that $w\notin N_{U_{\pi^*}}\big(v_p\big)$.

\noindent\textbf{Subcase 2.2~} $q\geq t+4$.
Since $d_q'\geq t+1$,
we have $d_q=d_q'+1\geq t+2$,
and so $d_{U_{\pi^*}}\big(v_q\big)=d_q-t\geq 2$.
Note that $N_{U_{\pi^*}}\big(v_q\big)$ contains exactly one vertex in $P_{U_{\pi^*}}\big(v_p,v_q\big)$,
and $q\geq t+4$.
Therefore,  there is vertex $w\in N_{U_{\pi^*}}\big(v_q\big)\setminus V(P_{U_{\pi^*}}\big(v_p,v_q)\big)$ such that
$w\notin N_{U_{\pi^*}}\big(v_p\big)\cup \big\{v_p\big\}$.

For both subcases above,
$w$ is a surprising vertex of $G_{\pi}$,
and hence $\Theta(G_{\pi},\alpha)<\Theta(G_{\pi'},\alpha)$ by Lemma \ref{54l}.
\end{pf}

\begin{thm}\label{57t}
If $c=2$,  $0\leq t\leq n-5$, and $\pi \lhd^* \pi'$ except for
$$\left((n-1)^{(t)},t+k+3,t+3,(t+2)^{(n-t-k-2)},(t+1)^{(k)}\right)
\lhd^*\left((n-1)^{(t)},t+k+4,(t+2)^{(n-t-k-1)},(t+1)^{(k)}\right),$$
where $1\leq k\leq n-t-5$,
then $\Theta(G_{\pi},\alpha)<\Theta(G_{\pi'},\alpha)$.
\end{thm}

\begin{pf}
According to Proposition \ref{41p},
we consider the following four cases.

\noindent\textbf{Case 1.~} $\pi$ satisfies $(4.2.1)$.
In this case,
$\pi=\left((n-1)^{(t)},t+4,(t+2)^{(n-t-1)}\right)$,
and $G\cong B_{\pi}^{(t)}=K_t\vee C(3,n-t-2)$ by Theorem \ref{46t}.
According to Theorem \ref{35t} and Lemma \ref{31l},
we may assume that
$N_H(v_{t+1})=\{w_{1},w_{2},w_{3},w_{4}\}$,
$\{w_{1},w_{2}\}\subseteq V(C_3)$, and
$\{w_{3},w_{4}\}\subseteq V(C_{n-t-2})$, and we may suppose that $f_{G_{\pi}}(w_{1})=\max\big\{f_{G_{\pi}}(w_{i})~|~1\leq i\leq 4\big\}.$

\noindent\textbf{Subcase 1.1~} $\pi'=\left((n-1)^{(t)},
t+5,(t+2)^{(n-t-2)},t+1\right)$.
Thus, $n-t\geq 6$,
and then there exists a vertex $x\in N_{C(3,n-t-2)}(w_{4})\setminus \big\{v_{t+1},w_{3}\big\}$.
Let $B_1=C(3,n-t-2)+v_{t+1}x-w_{4}x$ and $G_1=K_t\vee B_1$.
Clearly, $G_1\in \mathscr{B}(\pi',t)$.
Since $f_{G_{\pi}}(v_{t+1})>f_{G_{\pi}}(w_4)$ by Lemma \ref{31l},
then by Corollary \ref{23c} and the choice of $G_{\pi'}$,
we have $\Theta(G_{\pi},\alpha)<\Theta(G_1,\alpha)\leq\Theta(G_{\pi'},\alpha)$.

\noindent\textbf{Subcase 1.2~} $\pi'=\left((n-1)^{(t)},
t+4,t+3,(t+2)^{(n-t-2)},t+1\right)$.
Note that $n-t\geq 5$,
and there exists a vertex $x\in N_{C(3,n-t-2)}(w_{3})\setminus \{v_{t+1}\}$.
Let $B_2=C(3,n-t-2)+w_{1}w_{3}-xw_{3}$ and $G_2=K_t\vee B_1$.
Then, $G_2\in \mathscr{B}(\pi',t)$
and $f_{G_{\pi}}(w_{1})\geq f_{G_{\pi}}(x)$ by Theorem \ref{35t}.
Now, it follows from Corollary \ref{23c} and the choice of $G_{\pi'}$,
we have $\Theta(G_{\pi},\alpha)<\Theta(G_2,\alpha)\leq\Theta(G_{\pi'},\alpha)$.

\noindent\textbf{Case 2.~} $\pi$ satisfies $(4.2.2)$.
In this case,
$\pi=\left((n-1)^{(t)},(t+3)^{(2)},(t+2)^{(n-t-2)}\right)$,
and $G\cong B^{(t)}_{\pi}=K_t\vee H$, where $H=\theta(n-t-2,1,2)$ by Theorem \ref{46t}.
By Theorem \ref{35t} and Lemma \ref{31l},
we may assume that
$N_{H}(v_{t+1})=\{v_{t+2},w_{1},w_{2}\}$, $d_{H}(v_{t+2})=3$ and
$\big\{v_{t+1},v_{t+2},w_{1}\big\}$ induce a triangle in $H$.

\noindent\textbf{Subcase 2.1~} $\pi'=\left((n-1)^{(t)},
t+4,(t+2)^{(n-t-1)}\right)$.
Note that $n-t\geq 5$,
and then there exists a vertex $x\in N_{H}(v_{t+2})\setminus \big\{v_{t+1},w_{1},w_{2}\big\}$.
Let $B_3=H+v_{t+1}x-v_{t+2}x$ and $G_3=K_t\vee B_3$.
Clearly, $G_3\in \mathscr{B}(\pi',t)$,
and $f_{G_{\pi}}\big(v_{t+1}\big)\geq f_{G_{\pi}}\big(v_{t+2}\big)$ by Theorem \ref{35t}.

\noindent\textbf{Subcase 2.2~} $\pi'=\left((n-1)^{(t)},
(t+3)^{(3)},(t+2)^{(n-t-4)},t+1\right)$.
Note that $n-t\geq 5$,
and then there exists a vertex $x\in N_{H}(v_{t+2})\setminus \{v_{t+1},w_{1},w_{2}\}$.
Let $B_4=H+v_{t+2}w_{2}-v_{t+2}x$ and $G_4=K_t\vee B_4$.
It can be seen that $G_4\in \mathscr{B}(\pi',t)$,
and $f_{G_{\pi}}(w_{2})\geq f_{G_{\pi}}(x)$ by Theorem \ref{35t}.

\noindent\textbf{Subcase 2.3~} $\pi'=\left((n-1)^{(t)},
t+4,t+3,(t+2)^{(n-t-3)},t+1\right)$.
Note that $n-t\geq 5$,
and there exists a vertex $x\in N_{H}(w_{2})\setminus \{v_{t+1},v_{t+2}\}$.
Let $B_5=H+xv_{t+1}-xw_{2}$ and $G_5=K_t\vee B_5$.
Clearly, $G_5\in \mathscr{B}(\pi',t)$,
and $f_{G_{\pi}}(v_{t+1})>f_{G_{\pi}}(w_{2})$ by Theorem \ref{35t}.

For Subcases 2.1-2.3,
we always have $\Theta(G_{\pi},\alpha)<\Theta(G_{\pi'},\alpha)$
by Corollary \ref{23c} and the choice of $G_{\pi'}$, and so  the result follows.

\noindent\textbf{Case 3.~} $\pi$ satisfies $(4.2.3)$.
In this case,
$\pi=\left((n-1)^{(t)},t+k+4,(t+2)^{(n-t-k-1)},(t+1)^{(k)}\right)$,
where $1\leq k\leq n-t-5$.
By Theorem \ref{46t},
$G\cong B_{\pi}^{(t)}=K_t\vee B^\star$,
where $B^\star$ is the bicyclic graph with $n-t$ vertices obtained from $C(3,3)$ by attaching $d_{t+1}-4-t$ paths of almost equal lengths to the maximum degree vertex of $C(3,3)$.
By Theorems \ref{35t} and \ref{37t},
we may assume that
$N_H(v_{t+1})=\big\{v_{t+2},v_{t+3},v_{t+4},v_{t+5}\big\}$,
$v_{t+2}v_{t+3}\in E\big(C(3,3)\big)$ and $v_{t+4}v_{t+5}\in E\big(C(3,3)\big)$.

\noindent\textbf{Subcase 3.1~} $\pi'=\left((n-1)^{(t)},
t+k+4,t+3,(t+2)^{(n-t-k-3)},(t+1)^{(k+1)}\right)$.
It is easily checked  that $v_{t+5}\in N_{B^\star}(v_{t+4})\setminus \big(N_{B^\star}(v_{t+2})\cup \{v_{t+2}\}\big)$ and $f_{G_{\pi}}(v_{t+2})\geq f_{G_{\pi}}(v_{t+4})$ by Theorem \ref{35t}.
Let $B_6=B^\star+v_{t+5}v_{t+2}-v_{t+5}v_{t+4}$ and $G_6=K_t\vee B_6$.
Clearly, $G_6\in \mathscr{B}(\pi',t)$,
and then by Corollary \ref{23c} and the choice of $G_{\pi'}$, we get $\Theta(G_{\pi},\alpha)<\Theta(G_6,\alpha)\leq\Theta(G_{\pi'},\alpha)$.

\noindent\textbf{Subcase 3.2~} $\pi'=\left((n-1)^{(t)},
t+k+5,(t+2)^{(n-t-k-2)},(t+1)^{(k+1)}\right)$.
Thus,
$n-t-k-2\geq 4$,
and then $n-t-k-1\geq 5$.
Combining this with the structure of $B^\star$,
there exists a vertex $x\notin V\big(C(3,3)\big)$ with $d_{B^\star}(x)=2$,
and further there exists a vertex $y$ such that $y\in N_{B^\star}(x)$
and $y\notin V\big(P_{B^\star}(x,v_{t+1})\big)$.
Let $B_7=B^\star+v_{t+1}y-xy$ and $G_7=K_t\vee B_7$.
Obviously, $G_1\in \mathscr{B}(\pi',t)$
and $f_{G_{\pi}}(v_{t+1})\geq f_{G_{\pi}}(x)$ by Theorem \ref{35t}.
By Corollary \ref{23c} and the choice of $G_{\pi'}$, we get $\Theta(G_{\pi},\alpha)<\Theta(G_7,\alpha)\leq\Theta(G_{\pi'},\alpha)$.

\noindent\textbf{Case 4.~} $\pi$ satisfies $(4.2.4)$.
In this case, we have $G\cong B_\pi^{(t)}=K_t\vee B_{\pi^*}$ by Theorem \ref{46t} ,
where $B_{\pi^*}$ is the unique  BFS-bicyclic graph with $\pi^*$ as its degree sequence and $\mathcal{B}(B_{\pi^*})\cong \theta(2,1,2)$.
Following the notations in Remark \ref{53r},
we may suppose that $V\big(\theta(2,1,2)\big)=\big\{v_{t+1},v_{t+2},v_{t+3},v_{t+4}\big\}$, $d_{\theta(2,1,2)}\big(v_{t+1}\big)=d_{\theta(2,1,2)}\big(v_{t+2}\big)=3$,
$d_p'=d_p+1$,
$d_q=d_q'+1$, and
$d_i=d_i'$ for $i\in \big\{t+1,t+2,\ldots,n\big\}\setminus \big\{p,q\big\}$,
where $t+1\leq p<q\leq n$.

\noindent\textbf{Subcase 4.1~} $q\geq t+5$.
Since $d_q'\geq t+1$,
then $d_q=d_q'+1\geq t+2$,
and $d_{B_{\pi^*}}(v_q)=d_q-t\geq 2$.
Note that $N_{B_{\pi^*}}(v_q)$ contains exactly one vertex in $P_{B_{\pi^*}}(v_p,v_q)$,
and $q\geq t+5$.
Hence there is vertex $w\in N_{B_{\pi^*}}\big(v_q\big)\setminus V\big(P_{B_{\pi^*}}(v_p,v_q)\big)$ such that
$w\notin N_{B_{\pi^*}}\big(v_p\big)\cup \big\{v_p\big\}$.

\noindent\textbf{Subcase 4.2~} $q=t+3$ or $q=t+4$.
Notice that $d_{t+3}'\geq d_{t+4}'\geq t+2$ by $(4.2.4)$, as $\pi'$ is also a non-increasing degree sequence of $t$-cone bicyclic graph.
Then, $d_q=d_q'+1\geq t+3$.
Bearing in mind that $d_{B_{\pi^*}}(v_q)=d_q-t\geq 3$,
we have $v_p\in V\big(\theta(2,1,2)\big)$ and $N_{B_{\pi^*}}\big(v_q\big)\setminus V\big(\theta(2,1,2)\big)\neq \emptyset$.
So there exists a vertex $w\in N_{B_{\pi^*}}\big(v_q\big)\setminus V\big(\theta(2,1,2)\big)$ such that $w\notin N_{B_{\pi^*}}\big(v_p\big)\cup \big\{v_p\big\}$.

\noindent\textbf{Subcase 4.3~} $q=t+2$.
Then $p=t+1$.
According to $(4.2.4)$,
we have $d_{t+2}-t\geq 3$.
If $d_{t+2}=t+3$,
then $d'_{t+2}=t+2$,
and so it follows from  $(4.2.3)$ and  $(4.2.4)$ that
$$\pi'=\left((n-1)^{(t)},t+k+4,(t+2)^{(n-t-k-1)},(t+1)^{(k)}\right),\hspace{10pt}\text{and}$$
$$\pi=\left((n-1)^{(t)},t+k+3,t+3,(t+2)^{(n-t-k-2)},(t+1)^{(k)}\right),$$
where $1\leq k\leq n-t-5$,
contrary with  the condition.
Otherwise,  $d_{t+2}\geq t+4$, and hence   there exists a vertex $w\in N_{B_{\pi^*}}(v_{t+2})\setminus \left(N_{B_{\pi^*}}(v_{t+1})\cup \big\{v_{t+1}\big\}\right)$.

For all  subcases above,
$w$ is a surprising vertex of $G_\pi$,
and by Lemma \ref{54l} we have  $\Theta(G_\pi,\alpha)<\Theta(G_{\pi'},\alpha)$, as desired.
\end{pf}

At the end of this section,
we would like to point out that in Remark 3.2.2 of \cite{Liu1},
the condition $\pi \lhd^* \pi'$ should be added.


\section{Concluding Remark}

As a supplement of \cite{Luo}, by
taking $t=1$ and $\alpha\in \big\{0,~1\big\}$ in Theorems \ref{42t} and \ref{44t}, we have the  following Corollary \ref{61c}.
Besides, the following Corollary \ref{62c} can be obtained from Theorems \ref{55t} and \ref{56t}
by taking $t=1$ and $\alpha\in \big\{0,~1\big\}$,
which are the main results of \cite{Luo}.

\begin{cor}\label{61c}
Let $\pi$ be a non-increasing degree sequence of single-cone $c$-cyclic graph.

$(i)$ If $c=0$, then
$T_{\pi}^{(1)}$  is the unique $\rho$-maximal and the unique  $\mu$-maximal  graph in the class of single-cone trees  with $\pi$ as its degree sequence.

$(ii)$ If $c=1$, then
$U_{\pi}^{(1)}$  is the unique $\rho$-maximal and the unique  $\mu$-maximal  graph in the class of single-cone unicyclic graphs   with $\pi$ as its degree sequence.
\end{cor}

\begin{cor}\label{62c}{\em\cite{Luo}}
Let $G$ and $G'$ be the $\rho$-maximal $\big($resp., $\mu$-maximal$\big)$ graphs in the classes of single-cone $c$-cyclic graphs  with degree sequences $\pi$ and $\pi'$, respectively.
If $\pi\lhd \pi'$ and $c\in \big\{0,1\big\}$, then  $\rho(G)<\rho(G')$ $\big($resp.,~~$\mu(G)<\mu(G')\big)$.
\end{cor}

Moreover, by taking $t=0$ and $\alpha\in \big\{0,~1\big\}$ in Theorems \ref{42t}, \ref{44t}, \ref{46t}, \ref{55t} and \ref{56t}, we can deduce  the main results of \cite{Biyikoglu,Zhang0,FB2010,Liu2,Zhang1,Huang,LiuY2010} as follows:

\begin{cor}\label{63c}{\em\cite{Biyikoglu,Zhang0,FB2010,Zhang1,Huang,LiuY2010}}
Let $\pi$ be a non-increasing degree sequence of tree
$($resp. unicyclic graph, bicyclic graph$)$.
Then $T_{\pi}^{(0)}$ $\big($resp., $U^{(0)}_{\pi},\,B^{(0)}_{\pi}\big)$ is the unique $\rho$-maximal and the unique $\mu$-maximal graph in the class of trees
$\big($resp., unicyclic graphs, bicyclic graphs$\big)$ with $\pi$ as its degree sequence.
\end{cor}

\begin{cor}\label{64c}{\em\cite{Biyikoglu,Zhang0,Liu2,Zhang1}}
Let $G$ and $G'$ be the $\rho$-maximal $\big($resp., $\mu$-maximal$\big)$  graphs in the classes of $c$-cyclic graphs  with degree sequences $\pi$ and $\pi'$, respectively.
If  $\pi\lhd \pi'$ and $c\in \big\{0,1\big\}$, then  $\rho(G)<\rho(G')$ $\big($resp., $\mu(G)<\mu(G')\big)$.
\end{cor}

\end{document}